\RequirePackage{fix-cm}
\documentclass[smallextended]{svjour3}       
\smartqed  
\usepackage{graphicx}
\usepackage{amssymb}
\usepackage{latexsym}
\usepackage[overload]{empheq}

\usepackage{amsmath,amstext,amsthm,amssymb,bbm,amsbsy}
\usepackage{mathtools}
\usepackage{tabularx}
\usepackage{cite}
\usepackage[ruled,linesnumbered]{algorithm2e}
\usepackage{appendix}
\usepackage{multirow}

\usepackage[backref=page]{hyperref}
\renewcommand*{\backref}[1]{}
\renewcommand*{\backrefalt}[4]{%
    \ifcase #1 (Not cited.)%
    \or        (Cited on page~#2.)%
    \else      (Cited on pages~#2.)%
    \fi}
\newcommand{\be}{\begin{equation}}
\newcommand{\ee}{\end{equation}}
\newcommand{\benn}{\begin{equation*}}
\newcommand{\eenn}{\end{equation*}}

\DeclareMathOperator*{\argmin}{arg\,min}
\DeclareMathOperator*{\argmax}{arg\,max}

\newcommand{\bx}{\boldsymbol{x}}
\newcommand{\bu}{\boldsymbol{u}}
\newcommand{\bp}{\boldsymbol{p}}
\newcommand{\ba}{\boldsymbol{a}}
\newcommand{\bq}{\boldsymbol{q}}
\newcommand{\bv}{\boldsymbol{v}}
\newcommand{\by}{\boldsymbol{y}}
\newcommand{\bw}{\boldsymbol{w}}
\newcommand{\R}{\mathbb{R}}
\newcommand{\Rn}{\R^n}
\newcommand{\N}{\mathbb{N}}
\newcommand{\J}{\Psi}
\newcommand{\V}{V}
\newcommand{\bzero}{\mathbf{0}}
\newcommand{\symn}{S^{n}} 
\newcommand{\symposn}{S^{n}_{> 0}} 
\newcommand{\symposl}{S^{l}_{> 0}} 
\newcommand{\symnonnegn}{S^{n}_{\geq 0}} 
\newcommand{\matnn}{\R^{n\times n}}
\newcommand{\matnl}{\R^{n\times l}}
\newcommand{\On}{O_{n}}
\newcommand{\Ol}{O_{l}}

\newcommand{\Mxx}{Q} 
\newcommand{\Muu}{R} 
\newcommand{\Mxu}{S} 
\newcommand{\Gxx}{G} 
\newcommand{\Gx}{\ba} 
\newcommand{\Gc}{b} 
\newcommand{\Cpp}{C_{pp}} 
\newcommand{\Cxp}{C_{xp}} 
\newcommand{\Cxx}{C_{xx}} 
\newcommand{\Af}{A} 
\newcommand{\Bf}{B} 
\newcommand{\Sxx}{P} 
\newcommand{\Sx}{\bq} 
\newcommand{\Sc}{r} 
\newcommand{\constraint}{\mathcal{C}} 

\newcommand{\Soc}{\V}
\newcommand{\Voc}{\tilde{\V}}
\newcommand{\VNN}{V_{NN}}
\newcommand{\xzero}{\bx_0} 
\newcommand{\tzero}{t_0} 
\newcommand{\varx}{\bx} 
\newcommand{\vart}{t} 

\newcommand{\unn}{\bu_{NN}}
\newcommand{\optu}{\bu^*}
\newcommand{\optx}{\bx^*}
\newcommand{\opui}{\bv_i^*}

\newcommand{\opxi}{\by_i^*}

\begin{document}

\title{Neural network architectures using min-plus algebra for solving certain high dimensional optimal control problems and Hamilton-Jacobi PDEs\thanks{Research supported by AFOSR MURI FA9550-20-1-0358.
Authors' names are given in last/family name alphabetical order.}
}

\titlerunning{Neural networks for certain high dimensional optimal control problems}        

\author{J\'er\^ome Darbon         \and
        Peter M. Dower 
        \and Tingwei Meng
}

\institute{J\'er\^ome Darbon \at
              Division of Applied Mathematics, Brown University \\
              \email{jerome\_darbon@brown.edu}           
           \and
           Peter M. Dower \at
              Department of Electrical and Electronic Engineering, The University of Melbourne\\
              \email{pdower@unimelb.edu.au}   
        \and
        Tingwei Meng
        \at
      Department of Mathematics, UCLA \\
      \email{tingwei@math.ucla.edu}
}

\date{Received: date / Accepted: date}

\maketitle

\begin{abstract}
Solving high dimensional optimal control problems and corresponding Hamilton-Jacobi PDEs are important but challenging problems in control engineering. In this paper, we propose two abstract neural network architectures which are respectively used to compute the value function and the optimal control for certain class of high dimensional optimal control problems. We provide the mathematical analysis for the two abstract architectures. We also show several numerical results computed using the deep neural network implementations of these abstract architectures.
A preliminary implementation of our proposed neural network architecture on FPGAs shows promising speed up compared to CPUs.
This work paves the way to leverage efficient dedicated hardware designed for neural networks to solve high dimensional optimal control problems and Hamilton-Jacobi PDEs.
\keywords{Optimal control \and Hamilton--Jacobi partial differential equations \and Neural networks \and Grid-free numerical methods}
\end{abstract}

\section{Introduction}
Optimal control problems are an important class of optimization problems that find many applications in engineering, such as trajectory planning~\cite{Coupechoux2019Optimal,Rucco2018Optimal,Hofer2016Application,Delahaye2014Mathematical,Parzani2017HJB,Lee2021Hopf}, robot manipulator control~\cite{lewis2004robot,Jin2018Robot,Kim2000intelligent,Lin1998optimal,Chen2017Reachability} and humanoid robot control~\cite{Khoury2013Optimal,Feng2014Optimization,kuindersma2016optimization,Fujiwara2007optimal,fallon2015architecture,denk2001synthesis}. 
Under some assumptions, the optimal control problems are related to backward Hamilton--Jacobi (HJ) PDEs of the form
\begin{equation} \label{eqt: HJ}
\begin{dcases} 
-\frac{\partial \V}{\partial t}(t,\bx)+H(t,\bx, \nabla_{\bx}\V(t,\bx)) = 0 & \bx\in\mathbb{R}^{n}, t\in(0,T),\\
\V(T,\bx)=\J(\bx) & \bx\in\mathbb{R}^{n},
\end{dcases}
\end{equation}
where the function $H\colon [0,T]\times\Rn\times  \Rn\ni (t,\bx,\bp)\mapsto H(t,\bx,\bp)\in\R$ is called the Hamiltonian, which is convex with respect to $\bp$, and the continuous function $\J\colon\Rn\to\R$ specifies the terminal cost.
This relation between optimal control problems and HJ PDEs has been widely studied in the literature (see~\cite{Bardi1997Optimal} for instance). 
The value function of an optimal control problem may be characterized as the unique
viscosity solution of the corresponding HJ PDE,
while the optimal feedback control in the optimal control problems is related to the spatial gradient which is the Fréchet derivative of the value function with respect to the state variable $\bx$.
Therefore, computing the viscosity solution of an HJ PDE and its spatial gradient is an important problem in control engineering. 

In many practical engineering problems, the dimensionality is often high. For instance, in robot manipulator control problems, there are multiple joints and end effectors in the manipulator.
To control and measure the movement of each joint or end effector, several variables such as velocities, angles or positions are included in the state variable in the optimal control problems. 
As a result, the dimension of the state space is usually greater than five in practice. 
However, when the dimension is greater than five, standard grid-based numerical algorithms such as ENO \cite{Osher1991High}, WENO \cite{Jiang2000Weighted}, and DG \cite{Hu1999Discontinuous} are infeasible to apply.
This infeasibility is due to the curse of dimensionality~\cite{bellman1961adaptive}, i.e., as the dimension grows, the number of grid points grows exponentially, and hence the memory requirement as well as the computational time also grow exponentially.
Therefore, solving optimal control problems and HJ PDEs in high dimensions efficiently is an important but challenging problem. 
In the literature, several methods are proposed to overcome the curse of dimensionality when solving high dimensional HJ PDEs and optimal control problems. These methods include, but are not limited to, max-plus methods \cite{akian2006max,akian2008max, dower2015maxconference,Fleming2000Max,gaubert2011curse,mceneaney2006max,McEneaney2007COD,mceneaney2008curse,mceneaney2009convergence}, 
optimization methods \cite{darbon2015convex,darbon2019decomposition,Darbon2016Algorithms,yegorov2017perspectives}, 
tensor decomposition techniques \cite{dolgov2019tensor,horowitz2014linear,todorov2009efficient}, 
sparse grids \cite{bokanowski2013adaptive,garcke2017suboptimal,kang2017mitigating}, 
polynomial approximation \cite{kalise2019robust,kalise2018polynomial}, 
model order reduction \cite{alla2017error,kunisch2004hjb},
dynamic programming and reinforcement learning \cite{alla2019efficient,bertsekas2019reinforcement} and neural networks \cite{bachouch2018deep,bansal2020deepreach, Djeridane2006Neural,jiang2016using, Han2018Solving, hure2018deep, hure2019some, lambrianides2019new, Niarchos2006Neural, reisinger2019rectified,royo2016recursive, Sirignano2018DGM,Li2020generating,darbon2020overcoming,Darbon2021Neural,nakamurazimmerer2021adaptive,NakamuraZimmerer2021QRnet,jin2020learning,JIN2020Sympnets}. 

Recently, neural networks have been a successful tool in solving scientific computing problems involving PDEs. The related works include but are not limited to~\cite{bachouch2018deep,bansal2020deepreach,beck2018solving, beck2019deep, beck2019machine, Berg2018Unified,chan2019machine, Cheng2006Fixed, Djeridane2006Neural, Dissanayake1994Neural,  dockhorn2019discussion, E2017Deep, Farimani2017Deep, Fujii2019Asymptotic, grohs2019deep, Han2018Solving, han2019solving, hsieh2018learning, hure2018deep, hure2019some, jianyu2003numerical, khoo2017solving, khoo2019solving, Lagaris1998ANN, Lagaris2000NN, lambrianides2019new, lee1990neural, lye2019deep, McFall2009ANN, Meade1994Numerical, Milligen1995NN, Niarchos2006Neural, pham2019neural, reisinger2019rectified, royo2016recursive, Rudd2014Constrained, Sirignano2018DGM, Tang2017Study, Tassa2007Least, weinan2018deep, Yadav2015Intro, yang2018physics, yang2019adversarial, long2017pde,long2019pde, meng2019composite, meng2019ppinn, pang2019fpinns, raissi2018deep,raissi2018forward,raissi2017physicsi,raissi2017physicsii,Raissi2019PINN, uchiyama1993solving, zhang2019learning, zhang2019quantifying,Li2020generating,albi2021gradient,kang2020neural}. Due to the success of neural networks, many new hardware designs have been proposed to 
efficiently (in terms of speed, latency, throughput or energy) implement neural networks.
For instance, Google designed the ``Tensor Processor Unit" \cite{googleTPU17} to accelerate inference using neural networks, and Intel developed new specific low-level instructions in their processors to accelerate machine learning applications~\cite{banerjeeEtal2019sfi}. Field programmable gate arrays (FPGAs) have been successfully used to implement neural networks for real-time applications, see e.g.,~\cite{farabet-suml-11,farabet-fpl-09,farabet.09.iccvw}. There are also efforts for proposing completely new silicon designs~\cite{chen2020classification, Hirjibehedin.20.nature}  and  efficient hardware designs for standard activation functions~\cite{kundu2019ktanh}.  In addition, new computing architectures specialized for implementing neural network start to be available: for instance Xilinx recently launched a new computing architecture called Versal AI to efficiently implement neural networks. Note that these trends follow what LeCun suggested in~\cite[Sec.~3]{lecun2019isscc}. 
This dedicated hardware can in-principle
be used for any algorithm that can be represented as a neural network architecture.

Realising this new
hardware and silicon designs requires new dedicated software to implement neural networks on the new platforms. There are some available software development kits to convert neural network codes in standard frameworks such as PyTorch, TensorFlow, ONNX and HALO to executable codes on the aforementioned dedicated hardware. As long as an algorithm can be expressed in the neural network languages, it is possible to accelerate it with these new hardware.
Therefore, an algorithm must be expressed as a neural network to leverage these new computational platforms.

In the literature, most neural network based algorithms regard the space of neural networks as a finite dimensional function space which approximates abstract functional spaces in the problems, and this approximation is guaranteed by the universal approximation theorems (see~\cite{Kidger2020Universal,Leshno1993Multilayer,Hornik1991Approximation,Rossi2005Functional,Chen1993Approximations,lu2019deeponet} for instance).
The output is given by a neural network whose parameters are trained using a problem-related optimization model. However, in general, there is no guarantee for the convergence of the neural networks, and hence the outputs are not guaranteed to solve the targeted problems.
There is another research direction which focuses more on the neural network architectures, and provides theoretical guarantees for certain architectures. Along this research direction,~\cite{E2017Proposal,E2019meanfield} proposed the connections between Resnet architectures and numerical solvers for ODEs, and~\cite{darbon2020overcoming,Darbon2021Neural} presented several neural network architectures which express representation formulas for solving certain HJ PDEs.

In this work, we enlarge the class of HJ PDEs and optimal control problems which are solvable using neural network architectures by considering representation formulas for certain HJ PDEs with state and time dependent Hamiltonians.
We design the neural network architectures such that they solve the optimal control problems and
HJ PDEs of interest, with the neural network parameters assigned directly from the problem data, without the need for a training process.
Without the training process, our neural network architectures are guaranteed to solve the optimal control problems and HJ PDEs.

\textbf{Contributions of this paper.}
We present neural network architectures which solve certain high dimensional optimal control problems and the corresponding HJ PDEs. 
We consider the Hamiltonians $H$ in~\eqref{eqt: HJ} which are quadratic with respect to $(\bx,\bp)$ with coefficients depending on $t$, and the initial data $\J$ which is the minimum of finitely many quadratics. There are numerical solvers and theoretical analysis in the literature for these problems using linear-quadratic regulator and min-plus algebra. In this work, we present the neural network architectures according to these theories.
Our contribution is three-fold.
\begin{itemize}
    \item First, our work paves the way to leverage efficient dedicated hardware designed for neural networks to solve high dimensional optimal control problems and HJ PDEs. We present the neural network architectures based on the solid algorithms and theories in the literature for solving these problems. With the neural network architectures, it is possible to obtain efficient implementations in practice by converting the neural network codes to executable codes on the dedicated hardware. This facilitates future real-time implementations for solving high dimensional practical optimal control problems.
    Our work is easily implemented in standard frameworks, and we provide our implementations using TensorFlow in \url{https://github.com/TingweiMeng/NN_HJ_minplus}.
    \item 
    Unlike most neural network algorithms in the literature, we provide theoretical guarantees to prove that our neural network architectures solve certain optimal control problems and HJ PDEs. These theoretical guarantees follow from the linear-quadratic control problems and min-plus algebra techniques in the theories of optimal controls and HJ PDEs. In this way, we show the correspondence between optimal control theories and certain neural network architectures. This correspondence also provides possibilities for new interpretations of certain neural network architectures from the optimal control perspective.
    \item 
    We present an FPGA implementation of our proposed neural network architecture which shows that promising speed-ups can be expected compared to implementation on CPUs.
\end{itemize}

\textbf{Organization of this paper.}
The mathematical background of optimal controls and min-plus algebra is given in Section~\ref{sec:background}. In Section~\ref{sec:NNandHJandctrl}, two abstract neural network architectures are presented, which solve the HJ PDEs and are used to compute the optimal controls in the optimal control problems, respectively. 
The first abstract architecture is shown in Section~\ref{subsec: nn1} and depicted in Fig.~\ref{fig: nn_Riccati}, which is a one-layer neural network architecture with abstract neurons. It solves the HJ PDEs and the optimal values in the corresponding optimal control problems.
The second abstract architecture is shown in Section~\ref{subsec: nn2} and depicted in Fig.~\ref{fig: nn_u}, which is a two-layer neural network architecture with abstract neurons. It can be used to compute the optimal controls in the optimal control problems. 
In Section~\ref{sec:admm_method}, we consider more general terminal conditions and propose a numerical algorithm that combines our proposed neural network architecture and ADMM to solve the corresponding HJ PDEs and optimal control problems.
The implementations of our proposed two abstract architectures and their numerical results are presented in Section~\ref{sec:implementation}. There are different ways to implement the abstract architectures. Among these implementations we show the one using the fourth order Runge-Kutta method for illustration, which gives the deep Resnet-type implementations depicted in Figs.~\ref{fig: nn_V_RK4} and~\ref{fig: nn_u_RK4}. 
The numerical solutions computed by the proposed neural network architectures and implementations for three optimal control problems are shown in Sections~\ref{subsec: test_const},~\ref{subsec: test_tdep} and~\ref{subsec: test_newton}, respectively.
Section~\ref{subsec:example_admm} shows one numerical result with a general terminal condition.
An implementation  of our proposed neural network on a FPGA is described in Section~\ref{subsec:FPGA} and it shows promising speed-ups compared to a CPU implementation.
Some conclusions are drawn in Section~\ref{sec:conclusion}.

\section{Mathematical background} \label{sec:background}
Throughout, we use $\matnl$ to denote the set of matrices with $n$ rows and $l$ columns with entries in $\mathbb{R}$, and use $\symn$ to denote the set of real-valued symmetric matrices in $\matnn$. Also, $\symposn$ denotes the set of symmetric positive definite matrices in $\matnn$, and $\symnonnegn$ denotes the set of symmetric positive semi-definite matrices in $\matnn$.
We denote the identity matrix in $\matnn$ by $I_n$, and the zero matrix in $\matnn$ by $\On$. Moreover, we use the bold character to denote a vector, and we use the capital character to denote a matrix, if not mentioned specifically. The $\ell^2$-norm and $\ell^1$-norm in $\R^n$ are denoted by $\|\cdot\|$ and $\|\cdot\|_1$, respectively.

\subsection{Optimal control and min-plus algebra}
First, we give a brief introduction to optimal control problems, HJ PDEs and their relation.
An optimal control problem is formulated as follows
\begin{equation} \label{eqt: optctrl}
    \V(t_0,\bx_0) \doteq \inf\left\{\int_{t_0}^T L(s, \bx(s), \bu(s)) ds + \J(\bx(T))\right\}
\end{equation}
subject to 
\begin{equation}\label{eqt: ode}
    \begin{dcases}
    \dot{\bx}(s) = f(s, \bx(s), \bu(s)) & s\in (t_0,T),\\
    \bx(t_0) = \bx_0,
    \end{dcases}
\end{equation}
where $T\in (0,+\infty)$ and $t_0\in [0,T]$ are scalars which denote the terminal time and initial time, $\bx_0$ is a vector in $\Rn$ which denotes the initial position, the trajectory $\bx\colon [t_0,T]\to\Rn$ is an absolutely continuous function solving the Cauchy problem~\eqref{eqt: ode} almost everywhere, and the control $\bu\colon [t_0,T]\to \R^l$ is a function in a function space such as $L^p(t_0,T; \R^l)$ or the space of measurable functions. 
In the optimal control problem~\eqref{eqt: optctrl}, the running cost is given by the function $L\colon [0,T]\times\Rn\times \R^l\to \R$, which is also called Lagrangian, while the terminal cost is given by the function $\J\colon \Rn\to \R$.
The optimal cost is denoted by $\V(t_0,\bx_0)$, which is a function of the initial time $t_0$ and the initial position $\bx_0$ in the Cauchy problem~\eqref{eqt: ode}.

Under suitable assumptions (see~\cite{Bardi1997Optimal} for instance), the value function $\V$ is a viscosity solution of the corresponding backward HJ PDE~\eqref{eqt: HJ}. Viscosity solutions are known to be equivalent to minimax solutions (also known as minimal selections), which are defined via the
associated characteristic inclusion, see~\cite{Subbotin1996Minimax,Cannarsa2004Semiconcave}.
In the corresponding HJ PDE, the Hamiltonian $H\colon [0,T]\times\Rn\times \Rn \to \R \cup\{+\infty\}$ is given by the function $f$ and the Lagrangian $L$ as follows
\begin{equation*}
    H(t,\bx,\bp) = \sup_{\bu\in\R^l}\{ -\langle f(t,\bx,\bu),\bp\rangle - L(t,\bx,\bu)\} \quad \forall t\in [0,T], \bx,\bp\in\Rn, 
\end{equation*}
and the terminal data $\J\colon \Rn\to \R$ is given by the terminal cost in the optimal control problem.
Given the solution $\V$ to the HJ PDE~\eqref{eqt: HJ}, the optimal control $\bu^*\colon [0,T]\to\R^l$ in the problem~\eqref{eqt: optctrl} is characterized by Pontryagin maximum principle~\cite{Bardi1997Optimal}, which states that at almost every $s\in [0,T]$, the optimal control $\bu^*(s)$ satisfies
\begin{equation} \label{eqt: pmp}
    \bu^*(s) \in \argmax_{\bu\in\R^l}\{ -\langle f(s,\bx^*(s),\bu),\bp\rangle - L(s,\bx^*(s),\bu)\}, 
\end{equation}
for each $\bp\in D_{\bx}^+\V(s,\bx^*(s)) \cup D_{\bx}^-\V(s,\bx^*(s))$,
where $\bx^*\colon [t_0,T]\to\Rn$ is the corresponding trajectory solved by~\eqref{eqt: ode} given the control $\bu^*$. Here, $D_{\bx}^+\V$ and $D_{\bx}^-\V$ denote the set of the spatial components of the superdifferential and subdifferential of $V$, respectively.
There are different sets of assumptions for the above relation~\eqref{eqt: pmp} to hold. 
For details of the assumptions, see~\cite{Zhou1990Maximum},~\cite[Section~III.3.4]{Bardi1997Optimal} and the references in~\cite[Section~III.6]{Bardi1997Optimal}.
A verification theorem can alternatively be used to check whether a control is optimal, if a solution to the HJ PDE is known to exist.

We consider the HJ PDE~\eqref{eqt: HJ} whose terminal data $\J$ is the minimum of several functions $\J_i\colon\Rn\to\R$, i.e., we assume
\begin{equation} \label{eqt: Ji}
    \J(\bx) = \min_{i\in \{1,\dots, m\}} \J_i(\bx)\quad \forall \, \bx\in\Rn.
\end{equation}
Denote by $\V_i\colon \Rn\times [0,+\infty)\to \R$ the viscosity solution to the corresponding backward HJ PDE with terminal data $\J_i$, which reads
\begin{equation} \label{eqt: HJeqt_Ji}
\begin{dcases} 
-\frac{\partial \V_i}{\partial t}(t,\bx)+H(t,\bx,\nabla_{\bx}\V_i(t,\bx)) = 0 & \bx\in\mathbb{R}^{n}, t\in(0,T),\\
\V_i(T,\bx)=\J_i(\bx) & \bx\in\mathbb{R}^{n}.
\end{dcases}
\end{equation}
If the HJ PDEs~\eqref{eqt: HJ} and~\eqref{eqt: HJeqt_Ji} are solved by the value function~\eqref{eqt: optctrl} with terminal costs $\J$ and $\J_i$, respectively, then the solution operator in the HJ PDE~\eqref{eqt: HJ} is linear with respect to the min plus algebra~\cite{mceneaney2006max}. 
From straightforward calculation using~\eqref{eqt: optctrl}, the value function
$\V$ can be written as the minimum of $\V_i$ as follows
\begin{equation*} 
    \V(t,\bx) = \min_{i\in\{1,\dots,m\}} \V_i(t,\bx), \quad\forall \, \bx\in\Rn, t\in[0,T].
\end{equation*}
Consequently, $V$ in this form also solves the HJ PDE~\eqref{eqt: HJ} with terminal data $\J$.

\subsection{Neural networks}
We give a brief introduction of neural networks and refer the reader to~\cite{aggarwal2018neural} for a full introduction. 
A neural network architecture defines a space of functions which approximates the solution space in the target problem. 
A general neural network is the composition of several functions whose inputs and outputs are called layers. 
Each layer contains several variables or quantities which are called neurons. 
The input and output of the neural network function are called the input layer and the output layer, and all the other layers are called hidden layers. 
Different types of neural networks have been proposed in the literature~\cite{aggarwal2018neural}.

\begin{figure}[htbp]
\centering
\includegraphics[width = 0.7\textwidth]{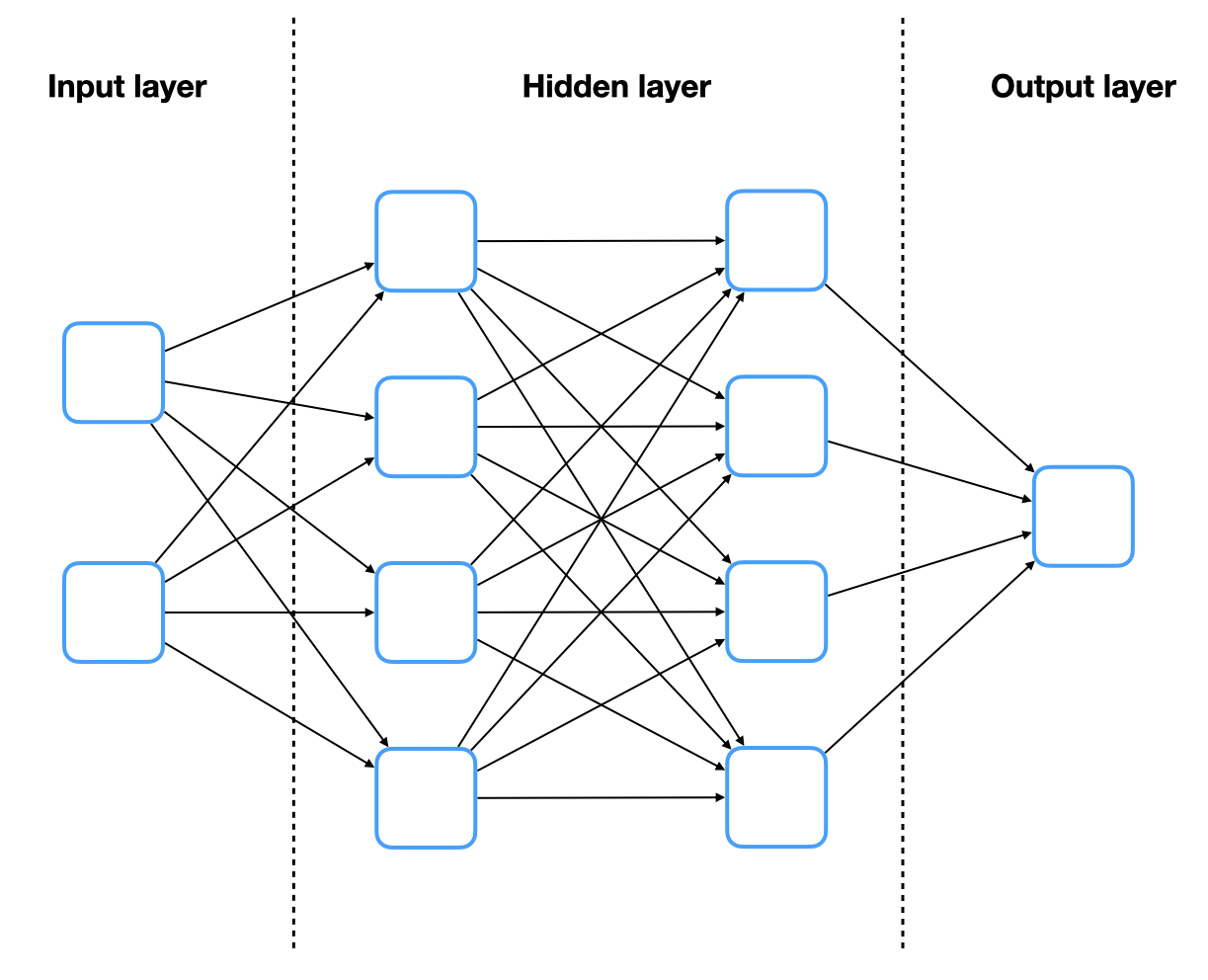}
\caption{
An illustration of a feedforward neural network architecture with two hidden layers.}
\label{fig:nn_general}
\end{figure}

A basic neural network architecture is called a feedforward neural network, whose hidden layer is the composition of an affine function and a non-linear function called activation function. 
An illustration of a feedforward neural network architecture with two hidden layers is shown in Fig.~\ref{fig:nn_general}, where each blue box corresponds to a neuron, and the line connecting the neurons illustrates the dependency between different neurons. To our knowledge, in the machine learning community, there is no standardised form to represent a neural network architecture as a diagram.

\begin{figure}[htbp]
\centering
\includegraphics[width = 0.5\textwidth]{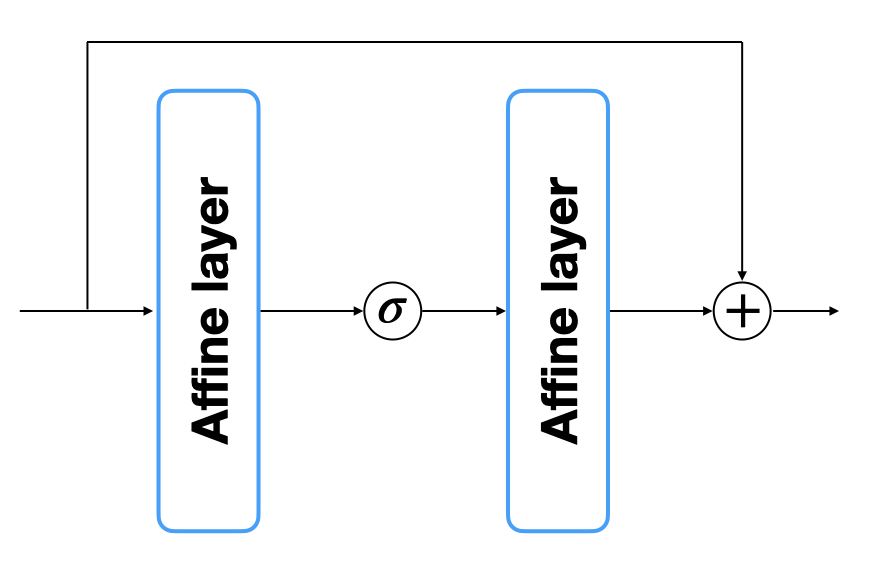}
\caption{
An illustration of a hidden layer in the Resnet architecture.}
\label{fig:nn_general_resnet}
\end{figure}

Another widely used neural network architecture is called residual neural network (Resnet)~\cite{He2016Deep,aggarwal2018neural}. 
A hidden layer in a Resnet involves more algebraic computations among compositions of affine functions and activation functions. 
An illustration of a hidden layer in a standard Resnet architecture with activation function $\sigma$ is shown in Fig.~\ref{fig:nn_general_resnet}.

\section{Neural network architectures for solving certain HJ PDEs and optimal control problems}\label{sec:NNandHJandctrl}
Fix a finite terminal time $T \in (0,+\infty)$, the following is assumed throughout:
\begin{itemize}
    \item[]  (A1)
    Let $\Af\in C([0,T];\matnn)$, $\Bf \in C([0,T];\matnl),\Mxu\in C([0,T];\matnl)$, $\Mxx\in C([0,T];\symposn)$ and $\Muu\in C([0,T];\symposl)$ be continuous functions. Let $\Gxx_i\in \symnonnegn$ be a constant matrix, $\Gx_i\in \Rn$ be a constant vector and $\Gc_i\in\R$ be a constant scalar for each $i\in\{1,\dots, m\}$.
\end{itemize}
We consider the optimal control problem~\eqref{eqt: optctrl} whose Lagrangian $L\colon [0,T]\times\Rn\times \R^l \to \R$ and the function $f\colon [0,T]\times\Rn\times \R^l \to \Rn$ are defined by
\begin{equation}\label{eqt:def_quad_L_f}
\begin{split}
    L(t,\bx,\bu) &= \frac{1}{2}\bx^T\Mxx(t)\bx + \frac{1}{2}\bu^T\Muu(t)\bu + \bx^T\Mxu(t)\bu,\\
    f(t,\bx,\bu) &= \Af(t)\bx + \Bf(t)\bu,
\end{split}
\end{equation}
for all $\bx\in\Rn$, $\bu\in \R^l$ and $t\in[0,T]$. Define the terminal cost $\J\colon\Rn\to\R$ by
\begin{equation} \label{eqt: H_J_quadratic_J}
\J(\bx) = \min_{i\in\{1,\dots, m\}} \left\{\frac{1}{2}\bx^T \Gxx_i\bx + \Gx_i^T\bx + \Gc_i\right\}    \quad\forall \bx\in\Rn.
\end{equation}
The corresponding optimal control problem is defined via the value function
\begin{equation} \label{eqt: optctrl_quadratic}
\begin{split}
    \Soc(\tzero, \xzero) = \inf_{(\bx,\bu)\in \constraint(\tzero,\xzero)}\Bigg\{\int_{\tzero}^T \Big(\frac{1}{2}\bx(s)^T\Mxx(s)\bx(s) 
    +  \frac{1}{2}\bu(s)^T\Muu(s)\bu(s) \quad \\
    +\, \bx(s)^T\Mxu(s)\bu(s)\Big) ds + \J(\bx(T))\Bigg\}
\end{split}
\end{equation}
where 
the constraint set $\constraint(\tzero,\xzero)$ is defined to be the set of $(\bx(\cdot),\bu(\cdot))\in L^2(\tzero,T; \Rn)\times L^2(\tzero,T;\R^l)$ which satisfies the following Cauchy problem
\begin{equation}\label{eqt: ode_quadratic}
    \begin{dcases}
    \dot{\bx}(s) = \Af(s)\bx(s)+\Bf(s)\bu(s) & s\in (\tzero,T),\\
    \bx(\tzero) = \xzero.
    \end{dcases}
\end{equation}
For the corresponding HJ PDE, we consider the following standard assumption~\cite{Yong1999Stochastic,Wang2014Deterministic}.
\begin{itemize}
    \item[] (A2)
    Assume $\Cpp\colon [0,T]\to\symnonnegn$, $, \Cxx\colon [0,T]\to\symposn$ and $\Cxp\colon [0,T]\to \matnn$ are three functions defined by
    \begin{equation} \label{eqt: H2}
        \begin{dcases}
        \Cpp(t) = \Bf(t)\Muu(t)^{-1}\Bf(t)^T,\\
        \Cxx(t) = \Mxx(t) - \Mxu(t)\Muu(t)^{-1}\Mxu(t)^T,\\
        \Cxp(t) = \Af(t) - \Bf(t)\Muu(t)^{-1}\Mxu(t)^T,
        \end{dcases}
    \end{equation}
    for all $t\in[0,T]$, where $\Af, \Bf, \Mxu,\Mxx,\Muu$ are the functions satisfying assumption (A1).
\end{itemize}
The Hamiltonian $H\colon [0,T]\times\Rn\times\Rn\to \R$ is defined by 
\begin{equation*} 
\begin{split}
    H(t,\bx,\bp) &= \frac{1}{2}\bp^T\Cpp(t)\bp - \frac{1}{2}\bx^T\Cxx(t)\bx - \bp^T \Cxp(t) \bx \quad \forall t\in[0,T], \bx,\bp\in\Rn.
\end{split}
\end{equation*}
The corresponding HJ PDE reads
\begin{equation} \label{eqt: HJ_quadratic_terminal}
\begin{dcases} 
-\frac{\partial \V(t,\bx)}{\partial t}+
H(t,\bx,\nabla_{\bx}\V(t,\bx)) = 0& \bx\in\mathbb{R}^{n}, t\in(0,T),\\
\V(T,\bx)=\J(\bx)=\min_{i\in\{1,\dots, m\}} \left\{\frac{1}{2}\bx^T \Gxx_i\bx + \Gx_i^T\bx + \Gc_i\right\}& \bx\in\mathbb{R}^{n}.
\end{dcases}
\end{equation}
This is the HJ PDE for a linear regulator problem with a min-of-quadratics terminal cost. Its solution can be represented via the solution to the Riccati equation (see~\cite{mceneaney2006max,dower2015maxconference,dower2016dynamic,Dower2016Game} for instance).

In the following two sections, we will present two abstract neural network architectures. The first one is shown in Section~\ref{subsec: nn1} which represents the viscosity solution to the HJ PDE~\eqref{eqt: HJ_quadratic_terminal}. The same neural network architecture also represents the value function in the optimal control problem~\eqref{eqt: optctrl_quadratic}.
The second abstract neural network architecture is shown in Section~\ref{subsec: nn2}, and it can be used to compute the optimal control in the optimal control problem~\eqref{eqt: optctrl_quadratic}.

\subsection{An abstract neural network architecture for solving the HJ PDE~\eqref{eqt: HJ_quadratic_terminal}} \label{subsec: nn1}

\begin{figure}[htbp]
\centering
\includegraphics[width = \textwidth]{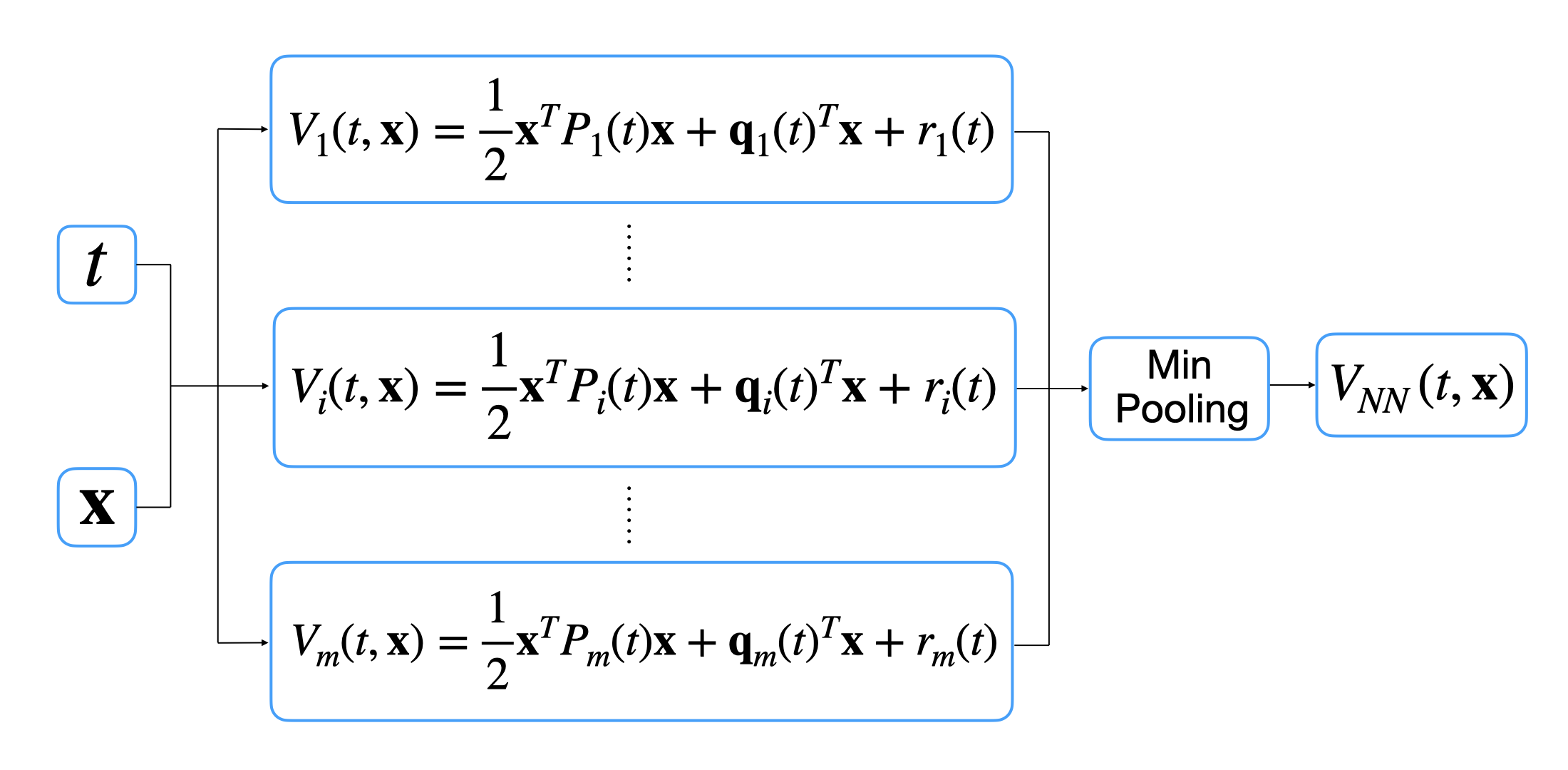}
\caption{
Illustration of the abstract neural network architecture defined by~\eqref{eqt: nn_S} that represents the viscosity solution to the HJ PDE~\eqref{eqt: HJ_quadratic_terminal} and the value function in the optimal control problem~\eqref{eqt: optctrl_quadratic}.}
\label{fig: nn_Riccati}
\end{figure}

We present an abstract neural network architecture which represents the viscosity solution to the HJ PDE~\eqref{eqt: HJ_quadratic_terminal}. The viscosity solution can be represented by a neural network $\VNN$ defined as follows
\begin{equation}\label{eqt: nn_S}
    \VNN(t,\bx) \doteq \min_{i\in\{1,\dots,m\}}\V_i(t,\bx), \quad \V_i(t,\bx) \doteq \frac{1}{2} \bx^T \Sxx_i(t)\bx + \Sx_i(t)^T\bx + \Sc_i(t),
\end{equation}
where the function $\Sxx_i\in C(0,T;\symn)$ solves the following Riccati final value problem (FVP)
\begin{equation} \label{eqt: odeP}
{\small
    \begin{dcases}
    \dot{\Sxx}_i(t) =  \Sxx_i(t)^T\Cpp(t)\Sxx_i(t) - \Sxx_i(t)^T\Cxp(t) - \Cxp(t)^T\Sxx_i(t) - \Cxx(t) &t\in(0,T)\\
    \Sxx_i(T) = \Gxx_i,
    \end{dcases}
    }
\end{equation}
the functions $\Sx_i\in C(0,T;\Rn)$ solves the following linear FVP
\begin{equation} \label{eqt: odeq}
    \begin{dcases}
    \dot{\Sx}_i(t) = \Sxx_i(t)^T\Cpp(t)\Sx_i(t) - \Cxp(t)^T\Sx_i(t) &t\in(0,T),\\
    \Sx_i(T) = \Gx_i,
    \end{dcases}
\end{equation}
and the function $\Sc_i\in C(0,T;\R)$ solves the following FVP
\begin{equation} \label{eqt: oder}
    \begin{dcases}
    \dot{\Sc}_i(t) = \frac{1}{2}\Sx_i(t)^T\Cpp(t)\Sx_i(t) &t\in(0,T),\\
    \Sc_i(T) = \Gc_i.
    \end{dcases}
\end{equation}
An illustration for the neural network architecture~\eqref{eqt: nn_S} is shown in Fig.~\ref{fig: nn_Riccati}. This is a one-layer abstract architecture with a min-pooling activation function. The $i$-th abstract neuron is given by the function $\V_i(t,\bx)$ in~\eqref{eqt: nn_S}. The architecture and neurons are called abstract since some ODE solvers for~\eqref{eqt: odeP},~\eqref{eqt: odeq} and~\eqref{eqt: oder} are further required in order to evaluate each neuron and the neural network architecture.
Later, in Section~\ref{sec:implementation}, we will provide a deep Resnet implementation (depicted in Fig.~\ref{fig: nn_V_RK4}) for this abstract architecture.

The following proposition shows that this neural network architecture~\eqref{eqt: nn_S} provides the viscosity solution to the HJ PDE~\eqref{eqt: HJ_quadratic_terminal} and the value function in the optimal control problem~\eqref{eqt: optctrl_quadratic}.

\begin{proposition}\label{prop: HJsol}
Assume (A1)-(A2) hold. Let $\VNN$ be the function defined by~\eqref{eqt: nn_S}. Then $\VNN$ is the unique viscosity solution to the HJ PDE~\eqref{eqt: HJ_quadratic_terminal}.
Moreover, $\VNN$ equals the value function $\Soc$ in the optimal control problem~\eqref{eqt: optctrl_quadratic}.
\end{proposition}

\begin{proof}
First, we apply~\cite[Prop.~2.2]{Wang2014Deterministic} to prove that the unique viscosity solution to the HJ PDE~\eqref{eqt: HJ_quadratic_terminal} is given by the value function for the optimal control problem~\eqref{eqt: optctrl_quadratic} under the assumptions (A1)-(A2). 
Most assumptions in~\cite[Prop.~2.2]{Wang2014Deterministic} are straightforward to check under this linear quadratic setting. Here, we only check the following three non-trivial assumptions: 
\begin{enumerate}
    \item There exists a positive constant $\rho>0$ such that $\Muu(t) - \rho I$ is positive semi-definite for all $t\in[0,T]$.
    \item There exist some positive constants $C_L>0$ and $\epsilon_0\in(0,1)$ such that
    $(1 - \epsilon_0)\bx^T\Mxx(t)\bx - \bx^T\Mxu(t)\Muu(t)^{-1}\Mxu(t)^T\bx \geq -C_L$ holds for all $t\in[0,T]$ and $\bx\in\Rn$.
    \item There exists a constant $C_0 >0$ such that $\J(\bx)\leq C_0(1+|\bx|^2)$ holds for all $\bx\in\Rn$.
\end{enumerate}
Now, we check the first assumption. 
We are going to prove a slightly stronger version. We prove that there exists $\rho>0$ such that the matrix $\Muu(t) - \rho I$ is positive definite for all $t\in[0,T]$.
Assume it does not hold. Then, there exist sequences $\{t_k\}\subset [0,T]$ and $\{\bx_k\}\subset \Rn$ such that $\bx_k^TR(t_k)\bx_k\leq \frac{1}{k}\|\bx_k\|^2$ holds for each $k\in\N$. After scaling, we assume $\|\bx_k\|=1$ holds for each $k\in\N$. By taking subsequences and still denoting the subsequences by $\{t_k\}$ and $\{\bx_k\}$, we obtain the convergence of $\{t_k\}$ and $\{\bx_k\}$, whose limits are denoted by $\bar{t}\in[0,T]$ and $\bar{\bx}\in\Rn$, respectively. Note that $\|\bx_k\|=1$ for each $k\in\N$ implies $\|\bar{\bx}\|=1$. Since $\Muu$ is a continuous function, we obtain
\begin{equation*}
    \bar{\bx}^TR(\bar{t})\bar{\bx} = \lim_{k\to\infty} \left(\bx_k^TR(t_k)\bx_k - \frac{1}{k}\|\bx_k\|^2\right) \leq 0,
\end{equation*}
which contradicts with the assumption that the matrix $R(\bar{t})$ is positive definite. Therefore, the first assumption holds.
A similar argument proves that there exists a positive constant $\epsilon$ such that $\Cxx(t) - \epsilon I$ is positive semi-definite for all $t\in[0,T]$. Also, the continuity of $\Mxx$ implies the existence of a uniform upper bound $C$ for 
$\|\Mxx(t)\|$ for all $t\in[0,T]$.
Let $\epsilon_0$ equal $\min\{\frac{\epsilon}{C},1\}$. Then, by~\eqref{eqt: H2}, we have
\begin{equation*}
\begin{split}
    (1 - \epsilon_0)\bx^T\Mxx(t)\bx - \bx^T\Mxu(t)\Muu(t)^{-1}\Mxu(t)^T\bx
    &= \bx^T \Cxx(x) \bx - \epsilon_0 \bx^T\Mxx(t)\bx \\
    &\geq\bx^T \Cxx(x) \bx - \epsilon_0 \|\Mxx(t)\|\|\bx\|^2\\
    &\geq \epsilon \|\bx\|^2 - \epsilon_0 C \|\bx\|^2 \geq 0.
\end{split}
\end{equation*}
As a result, the second assumption holds. The third assumption follows from a straightforward computation which reads
\begin{equation*}
\begin{split}
    \J(\bx) &= \min_{i\in\{1,\dots, m\}} \left\{\frac{1}{2}\bx^T \Gxx_i\bx + \Gx_i^T\bx + \Gc_i\right\}
    \leq \frac{1}{2}\bx^T \Gxx_1\bx + \Gx_1^T\bx + \Gc_1\\
    &\leq C_0(1+\|\bx\|^2),
\end{split}
\end{equation*}
for some positive constant $C_0$.
Therefore, all the assumptions in~\cite[Prop.~2.2]{Wang2014Deterministic} are satisfied. Then, by~\cite[Prop.~2.2]{Wang2014Deterministic}, the value function defined by~\eqref{eqt: optctrl_quadratic} is the unique viscosity solution to the HJ PDE~\eqref{eqt: HJ_quadratic_terminal}. (Note that if all the functions in (A1) and (A2) do not depend on time $t$, then this result follows from \cite{Bardi1997Bellman}.)
Therefore, it suffices to prove that $\VNN$ defined by~\eqref{eqt: nn_S} is the value function in the optimal control problem~\eqref{eqt: optctrl_quadratic}, i.e., it suffices to prove that $\VNN(\tzero,\xzero) =\Soc(\tzero, \xzero)$ holds for all $\xzero\in\Rn$ and $\tzero\in [0,T]$.

Define $\J_i\colon \Rn\to\R$ by 
\begin{equation}\label{eqt: defJi}
\J_i(\bx) \doteq \frac{1}{2}\bx^T \Gxx_i\bx + \Gx_i^T\bx + \Gc_i, \quad \forall\, \bx\in\Rn,
\end{equation}
for each $i\in\{1,\dots,m\}$.
Then, by~\eqref{eqt: H_J_quadratic_J}, the functions $\J_1,\dots, \J_m$ and $\J$ satisfy~\eqref{eqt: Ji}.
According to~\eqref{eqt: optctrl_quadratic},~\eqref{eqt: Ji}
and the min-plus linearity of the dynamic programming evolution operator, we have
\begin{equation}\label{eqt: minplus_pf31}
\begin{split}
    &\Soc(\tzero, \xzero) \\
    = &\inf_{(\bx,\bu)\in \constraint(\tzero,\xzero)}\Bigg\{\int_{\tzero}^T \Big(\frac{1}{2}\bx(s)^T\Mxx(s)\bx(s) +  \frac{1}{2}\bu(s)^T\Muu(s)\bu(s) \\
    &\quad\quad\quad\quad\quad\quad\quad\quad\quad\quad + \bx(s)^T\Mxu(s)\bu(s)\Big) ds + \min_{i\in\{1,\dots, m\}}\J_i(\bx(T))\Bigg\}\\
    = &\min_{i\in\{1,\dots, m\}}\Bigg\{\inf_{(\bx,\bu)\in \constraint(\tzero,\xzero)}\Bigg\{\int_{\tzero}^T \Big(\frac{1}{2}\bx(s)^T\Mxx(s)\bx(s) +  \frac{1}{2}\bu(s)^T\Muu(s)\bu(s) \\
    &\quad\quad\quad\quad\quad\quad\quad\quad\quad\quad\quad\quad\quad + \bx(s)^T\Mxu(s)\bu(s)\Big) ds + \J_i(\bx(T))\Bigg\}\Bigg\}.
\end{split}
\end{equation}
We define the value function in the last line of~\eqref{eqt: minplus_pf31} to be $\Voc_i(\tzero,\xzero)$, i.e., we define the function $\Voc_i\colon [0,T]\times \Rn \to\R$ by
\begin{equation} \label{eqt: optctrl_quadratic_i}
\begin{split}
    \Voc_i(\tzero, \xzero) \doteq \inf_{(\bx,\bu)\in \constraint(\tzero,\xzero)}\Bigg\{\int_{\tzero}^T \Big(\frac{1}{2}\bx(s)^T\Mxx(s)\bx(s) +  \frac{1}{2}\bu(s)^T\Muu(s)\bu(s) \quad \\
    + \,\bx(s)^T\Mxu(s)\bu(s)\Big) ds + \J_i(\bx(T))\Bigg\}.
\end{split}
\end{equation}

It is well-known that $\Voc_i$ can be solved by Riccati equation under the assumptions (A1)-(A2). One way to prove it is given as follows. First,~\cite{Bucy1967Global} shows the existence of the global solution to the Riccati FVP~\eqref{eqt: odeP} under the assumptions (A1)-(A2). Then, with a similar argument as in the proof of~\cite[Chap~6, Thm.~2.8]{Yong1999Stochastic}, the value function $\Voc_i$ is proved to satisfy 
\begin{equation}\label{eqt: Si_quadratic}
    \Voc_i(t,\bx) = \frac{1}{2} \bx^T \Sxx_i(t)\bx + \Sx_i(t)^T\bx + \Sc_i(t) = \V_i(t,\bx),
\end{equation}
for all $\bx\in\R^n$ and $t\in[0,T]$, where $\V_i$ is the function defined in~\eqref{eqt: nn_S}, and $\Sxx_i$, $\Sx_i$ and $\Sc_i$ satisfy~\eqref{eqt: odeP},~\eqref{eqt: odeq}, and~\eqref{eqt: oder}, respectively.

Therefore, combining~\eqref{eqt: minplus_pf31} and~\eqref{eqt: Si_quadratic}, we derive that $\VNN\equiv \Soc$ in $\Rn\times[0,T]$ and the conclusion follows.
\end{proof}

\def\knn{k_{NN}}

\subsection{An abstract neural network architecture for the optimal control problem~\eqref{eqt: optctrl_quadratic}} \label{subsec: nn2}

\begin{figure}[htbp]
\centering
\includegraphics[width = \textwidth]{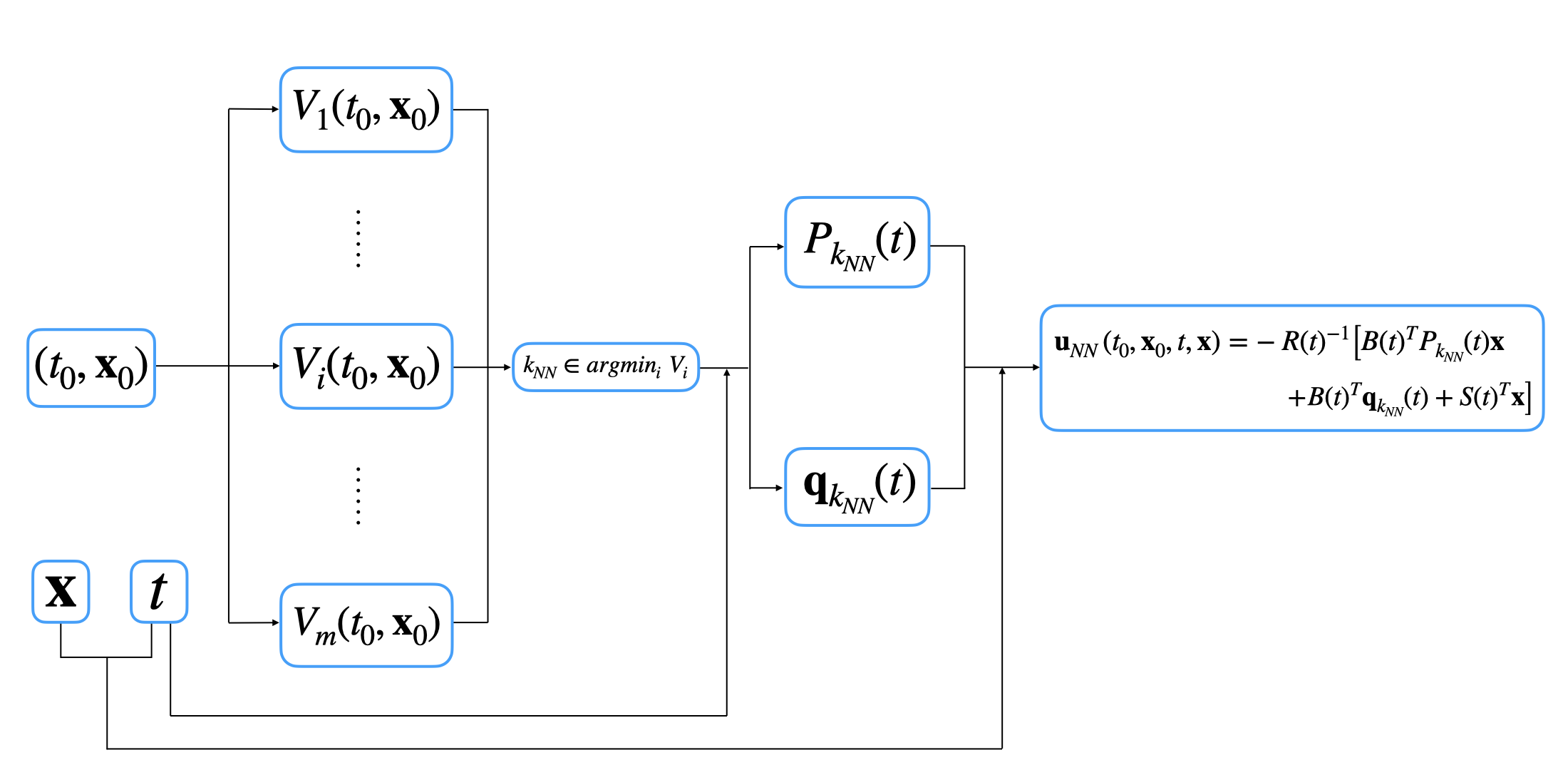}
\caption{
Illustration of the abstract neural network architecture defined by~\eqref{eqt: nn_u} that can be used to compute the optimal control in the optimal control problem~\eqref{eqt: optctrl_quadratic}.}
\label{fig: nn_u}
\end{figure}

We present an abstract neural network architecture for computing the optimal control in the optimal control problem~\eqref{eqt: optctrl_quadratic}. 
For any fixed index $j\in\{1,\dots, m\}$, define a function $\bu_j\colon [t_0,T]\times \Rn\to\R^l$ by
\begin{equation}\label{eqt:def_uj}
    \bu_j(\vart, \varx) \doteq -\Muu(\vart)^{-1}\left(\Bf(\vart)^T\Sxx_j(\vart)\varx + \Bf(\vart)^T\Sx_j(\vart) + \Mxu(\vart)^T\varx\right),
\end{equation}
for any $t\in [t_0,T]$ and $\bx\in\Rn$.
It is well-known that the function $\bu_j$ is the feedback control for the optimal control problem~\eqref{eqt: optctrl_quadratic_i} with index $i=j$. 
To solve the optimal control problem~\eqref{eqt: optctrl_quadratic}, we define the following function $\unn$ with specific selected index $\knn$:
\begin{equation}\label{eqt: nn_u}
    \unn(t_0,\bx_0, t,\bx) = \bu_{\knn(t_0,\bx_0)}(t,\bx),\quad \forall \bx_0,\bx\in\Rn, 0\leq t_0\leq t\leq T,
\end{equation}
where the index function $k_{NN}\colon [0,T]\times \Rn\to \{1,\dots,m\}$ is defined by
\begin{equation} \label{eqt: defk}
\begin{split}
    k_{NN}(t_0,\bx_0) &\in\argmin_{i\in\{1,\dots, m\}} \left\{\frac{1}{2} \xzero^T \Sxx_i(\tzero)\xzero + \Sx_i(\tzero)^T\xzero + \Sc_i(\tzero)\right\} \\
    &= \argmin_{i\in\{1,\dots, m\}} \V_i(\tzero,\xzero).
\end{split}
\end{equation}
When there is no ambiguity, we abuse the notation $k_{NN}(t_0,\bx_0)$ with $\knn$.
Recall that the functions $\Sxx_i$, $\Sx_i$, $\Sc_i$ and $\V_i$ are the functions defined in~\eqref{eqt: odeP},~\eqref{eqt: odeq}, \eqref{eqt: oder} and~\eqref{eqt: nn_S}, respectively.
If there are more than one minimizer in the optimization problem in~\eqref{eqt: defk}, we just select any of these minimizers and it will provide an optimal control. We will discuss more about this non-uniqueness later in Remark~\ref{rem:nonunique_control}.
Note that the function $\unn$ can be expressed using an abstract neural network architecture shown in Fig.~\ref{fig: nn_u}.
As it is discussed in Section~\ref{subsec: nn1}, the evaluation of the neurons $\{\V_i(\tzero,\xzero)\}_{i=1}^m$, $P_{\knn}(t)$ and $\bq_{\knn}(t)$ (where $\knn$ is the index defined in~\eqref{eqt: defk}) require further ODE solvers, and hence we call it an abstract architecture. An implementation of this abstract architecture using deep Resnet neural networks is provided in Section~\ref{sec:implementation} and depicted in Fig.~\ref{fig: nn_u_RK4}.

For completeness, the following proposition proves that the function $\unn$ defined in~\eqref{eqt: nn_u} computes the optimal control in the problem~\eqref{eqt: optctrl_quadratic}
if $(t,\bx)$ is on the optimal trajectory with initial time $t_0$ and initial position $\bx_0$.

\begin{proposition}\label{prop: ctrl}
Let $\bx_0\in\Rn$ and $t_0\in [0,T]$ be the initial position and initial time. Assume (A1)-(A2) hold. Let $\optu\colon [\tzero, T]\to\R^l$ be a feasible control, and $\optx\in C([\tzero,T]; \Rn)$ be the solution to the Cauchy problem~\eqref{eqt: ode_quadratic} with the control $\optu$.
Then, the function $\optu$ is an optimal control in the problem~\eqref{eqt: optctrl_quadratic} if and only if there exists an index $k$ in the set of minimizers of the optimization problem in~\eqref{eqt: defk}, 
such that there holds
    \begin{equation} \label{eqt: prop2_feedback}
    \optu(\vart)= \bu_k(\vart,\optx(\vart)) \quad \forall \vart\in[\tzero,T],
    \end{equation}
    where the function $\bu_k\colon [\tzero,T]\times\R^n \to \R^l$ is defined by~\eqref{eqt:def_uj} with the index $k$.
\end{proposition}
\begin{proof}
In this proof, we adopt the same notations as in the proof of Prop.~\ref{prop: HJsol}, but we abuse the notation $\VNN$ with $\V$ because they are shown to be equal in the proof of Prop.~\ref{prop: HJsol}. Let $\J_i$ be the function defined by~\eqref{eqt: defJi}. We consider the corresponding optimal control problem~\eqref{eqt: optctrl_quadratic_i}.
It is well known that the optimal control of the linear quadratic optimal control problem~\eqref{eqt: optctrl_quadratic_i} exists and is unique under the assumptions (A1)-(A2). We denote the optimal control by $\opui\colon [\tzero,T]\to\R^l$, and denote the corresponding trajectory by $\opxi\colon [\tzero,T]\to\Rn$. Moreover, the feedback form of the optimal control is given by 
\begin{equation}\label{eqt: eqt_vi_pf2}
\begin{split}
    \opui(\vart) &= 
    -\Muu(\vart)^{-1}(\Bf(\vart)^T\nabla_{\bx}\V_i(\vart,\opxi(\vart)) + \Mxu(\vart)^T\opxi(\vart))\\
     &= -\Muu(\vart)^{-1}(\Bf(\vart)^T\Sxx_i(\vart)\opxi(\vart) + \Bf(\vart)^T\Sx_i(\vart) + \Mxu(\vart)^T\opxi(\vart))\\
     &= \bu_i(t, \opxi(\vart)),
\end{split}
\end{equation}
for any $\vart\in[\tzero,T]$, where $\V_i$ is defined by~\eqref{eqt: Si_quadratic}. This is proved, for instance, by~\cite{Bucy1967Global} and a similar argument in the proof of~\cite[Chap~6, Thm.~2.8]{Yong1999Stochastic}.

Now, we prove the first implication of the proposition statement. Let $k$ be an index satisfying~\eqref{eqt: defk}, and the function 
$\bu_k\colon [\tzero,T]\times\R^n \to \R^l$ be defined by~\eqref{eqt:def_uj} with the index $k$. 
From~\eqref{eqt: eqt_vi_pf2}, we conclude that
the open loop control $t\mapsto \optu(\vart)=\bu_k(\vart,\optx(\vart))$ is the optimal control for the problem~\eqref{eqt: optctrl_quadratic_i} with the index $k$. Moreover, by~\eqref{eqt: minplus_pf31} and the definition of $k$, the optimal value $\V(t_0,\bx_0)$ in the problem~\eqref{eqt: optctrl_quadratic} equals the optimal value $\V_k(t_0,\bx_0)$ in the problem~\eqref{eqt: optctrl_quadratic_i} with the index $k$. Therefore, $t\mapsto \optu(\vart)$ is also an optimal control for the problem~\eqref{eqt: optctrl_quadratic}. 

Then, we prove the other implication. Assume $\optu$ is an optimal control in the problem~\eqref{eqt: optctrl_quadratic}, and $\optx$ is the corresponding optimal trajectory. Let $\tilde{k}$ be an index satisfying 
\begin{equation*}
    \tilde{k}\in \argmin_{i\in\{1,\dots,m\}} \J_i(\optx(T)).
\end{equation*}
Then, we have
\begin{equation}\label{eqt: pf_prop32_ineqtS}
\begin{split}
\V(t_0,\bx_0) &= \int_{\tzero}^T \Big(\frac{1}{2}\optx(\vart)^T\Mxx(\vart)\optx(\vart) + \frac{1}{2}\optu(\vart)^T\Muu(\vart)\optu(\vart) \\
&\quad\quad\quad\quad\quad\quad\quad+ \optx(\vart)^T\Mxu(\vart)\optu(\vart)\Big) d\vart + \J(\optx(T))\\
&= \int_{\tzero}^T \Big(\frac{1}{2}\optx(\vart)^T\Mxx(\vart)\optx(\vart) + \frac{1}{2}\optu(\vart)^T\Muu(\vart)\optu(\vart) \\
&\quad\quad\quad\quad\quad\quad\quad+ \optx(\vart)^T\Mxu(\vart)\optu(\vart)\Big) d\vart + \J_{\tilde{k}}(\optx(T))\\
&\geq \V_{\tilde{k}}(t_0,\bx_0)
\geq \V(t_0,\bx_0),
\end{split}
\end{equation}
where the first equality holds since $\optu$ is an optimal control of~\eqref{eqt: optctrl_quadratic} with the trajectory $\optx$, the second equality holds by definition of $\tilde{k}$, the first inequality holds since $\V_{\tilde{k}}(t_0,\bx_0)$ is the optimal value of the problem~\eqref{eqt: optctrl_quadratic_i} with the index $\tilde{k}$, and the last inequality holds by~\eqref{eqt: minplus_pf31}.
As a result, the two inequalities in~\eqref{eqt: pf_prop32_ineqtS} both become equalities, which implies that $\optu$ is the optimal control of the problem~\eqref{eqt: optctrl_quadratic_i} with the index $\tilde{k}$, and $\tilde{k}$ is a minimizer of the optimization problem in~\eqref{eqt: defk}.
Recall that the unique optimal control of the problem~\eqref{eqt: optctrl_quadratic_i} with the index $\tilde{k}$ satisfies the feedback form~\eqref{eqt: eqt_vi_pf2}, and hence we get~\eqref{eqt: prop2_feedback}.
\end{proof}

\begin{remark} \label{rem:nonunique_control}
Note that the existence of the optimal control $\optu$ is given by Prop.~\ref{prop: ctrl}, the existence of $\knn$ in~\eqref{eqt: defk}, and the existence of the optimal control with terminal cost $\J_{\knn}$.
However, such $\optu$ may not be unique, since there may be more than one minimizer in~\eqref{eqt: defk}. 
It can be seen in the above proposition any minimizer $k$ in~\eqref{eqt: defk} can define an optimal control
in the problem~\eqref{eqt: optctrl_quadratic}. As a result, if $\argmin_{i\in\{1,\dots, m\}} \V_i(\tzero, \xzero)$ is not a singleton, then there may be more than one optimal control in the problem~\eqref{eqt: optctrl_quadratic}. This non-uniqueness is possible since the optimal control problem is a non-convex optimization problem, where the terminal condition $\J$ is non-convex. 
For a fixed initial position $\bx_0$ and initial time $t_0$, one candidate of open loop optimal control is $\bv^*_{\knn}$ in~\eqref{eqt: eqt_vi_pf2} with index $i=\knn$.
If there are more than one minimizer in~\eqref{eqt: defk}, we select $\knn$ to be one minimizer, and then an optimal control $\optu$ is computed using $\unn$ in~\eqref{eqt: nn_u}.
\end{remark}

\begin{remark} \label{rem:compute_control}
We can compute the open-loop optimal control using our proposed neural network architecture $\unn$. 
For a fixed initial time $t_0$ and initial position $\bx_0$, we combine the function $\unn$ with the Cauchy problem~\eqref{eqt: ode_quadratic} and obtain the following Cauchy problem
\begin{equation} \label{eqt: ODE_optx}
    \begin{dcases}
    \dot{\bx}(s) = \Af(s)\bx(s)+\Bf(s)\unn(t_0, \bx_0, s,\bx(s)) & s\in (\tzero,T),\\
    \bx(\tzero) = \xzero.
    \end{dcases}
\end{equation}
By straightforward calculation using~\eqref{eqt: H2} and~\eqref{eqt: nn_u}, the differential equation above becomes
\begin{equation*}
{\small
\begin{split}
    \dot{\bx}(s) &= \Af(s)\bx(s)+\Bf(s)\unn(t_0, \bx_0, s,\bx(s))\\
    &= \Af(s)\bx(s)-\Bf(s)\Muu(s)^{-1}\left(\Bf(s)^T\Sxx_{\knn}(s)\varx(s) + \Bf(s)^T\Sx_{\knn}(s) + \Mxu(s)^T\varx(s)\right)\\
    &= \left(\Cxp(s) - \Cpp(s)\Sxx_{\knn}(s)\right)\bx(s) - \Cpp(s)\Sx_{\knn}(s).
\end{split}
}
\end{equation*}
The solution to this Cauchy problem is the optimal trajectory $\optx$. From the optimal trajectory, we obtain the open loop optimal control $\optu$ by $\optu = \unn(t_0, \bx_0, s,\optx(s))$. Our numerical results in Section~\ref{sec:implementation} are computed using this procedure.
\end{remark}

\begin{remark}
Note that the function $\unn$ also computes the feedback optimal control. Consider the function $(t,\bx)\mapsto \unn(t,\bx, t, \bx)$. Applying Prop.~\ref{prop: ctrl} to $\bx=\bx_0$ and $t=t_0$, we conclude that the optimal control at $t_0$ with initial time $t_0$ and position $\bx_0$ is $\optu(t_0) = \unn(t_0, \bx_0, t_0,\bx_0)$. Since the initial time $t_0\in[0,T]$ and initial position $\bx_0\in\Rn$ can be arbitrary, we conclude that function $(t,\bx)\mapsto \unn(t,\bx, t, \bx)$ is a feedback optimal control for problem~\eqref{eqt: optctrl_quadratic}.
\end{remark}

\def\indk{j}

\subsection{The extension to general terminal costs}\label{sec:admm_method}
Now, we consider more general terminal costs $\J$ of the form
\begin{equation} \label{eqt:J_min_Ji}
    \J(\bx) = \min_{i\in\{1,\dots,m\}} \J_i(\bx) \quad \forall \bx\in\Rn,
\end{equation}
where $\J_1,\dots, \J_m$ are some convex functions whose proximal points are numerically computable. Note that we no longer assume $\J_1,\dots, \J_m$ to be quadratics.
By min-plus linearity, the solution $\V$ is given by
\begin{equation}
    \V(\bx,t) = \min_{i\in\{1,\dots,m\}} \V_i(\bx,t) \quad\forall \bx\in\Rn, t\geq 0,
\end{equation}
where each $\V_i$ is the value function of the optimal control problem with terminal cost $\J_i$. Therefore, to solve this problem, we need to solve the $m$ subproblems. In the $i$-th subproblem, we solve the optimal control problem with terminal cost $\J_i$. If $\J_i$ is in the quadratic form, we apply the aforementioned method and neural network architecture to solve it. Otherwise, we apply the ADMM method~\cite{Glowinski2014Alternating,Boyd2011Distributed} to solve the $i$-th subproblem, whose $\indk$-th iteration includes the following three steps:
\begin{enumerate}
    \item Solve the optimal control problem
    \begin{equation} \label{eqt:admm_step1}
    \begin{split}
        \min_{(\bx,\bu)\in \constraint(\tzero,\xzero)}\Bigg\{\int_{\tzero}^T \frac{1}{2}\bx(s)^T\Mxx(s)\bx(s) 
        +  \frac{1}{2}\bu(s)^T\Muu(s)\bu(s) \quad \\
        +\, \bx(s)^T\Mxu(s)\bu(s) ds + \frac{\rho}{2}\|\bx(T) - \by^{\indk,i} + \bw^{\indk,i}\|^2\Bigg\}
    \end{split}
    \end{equation}
    where the constraint set $\constraint(\tzero,\xzero)$ contains the solutions to~\eqref{eqt: ode_quadratic}.
    Denote the optimal trajectory by $\bx^{\indk+1,i}(\cdot)$. Note that we add the superscript $i$ the emphasize that these minimizers are used to solve the $i$-th subproblem. 
    Since the terminal cost in~\eqref{eqt:admm_step1} is a quadratic function, this optimal control problem can be solved using Riccati equation, which can be represented by the neural network architecture depicted in Fig.~\ref{fig: nn_Riccati} with one neuron (i.e., $m=1$). The solution is given by $\frac{1}{2} \bx^T \Sxx_i(t)\bx + \Sx_i(t)^T\bx + \Sc_i(t)$, where $\Sxx_i(t)$, $\Sx_i(t)$ and $\Sc_i(t)$ solve the FVPs~\eqref{eqt: odeP},~\eqref{eqt: odeq} and~\eqref{eqt: oder} whose terminal conditions are given by 
    \begin{equation*}
    \Sxx_i(T) = \rho I_n, \quad \Sx_i(T) = \rho(\bw^{\indk,i} - \by^{\indk,i}), \quad \Sc_i(T) = \frac{\rho}{2}\|\bw^{\indk,i}-\by^{\indk,i}\|^2.    
    \end{equation*}
    \item Solve the following proximal point problem
    \begin{equation*}
        \min_{\by\in\Rn} \J_i(\by) + \frac{\rho}{2}\|\bx^{\indk+1,i}(T) - \by + \bw^{\indk,i}\|^2,
    \end{equation*}
    and denote the minimizer by $\by^{\indk+1,i}$. 
    \item Update $\bw$ by $\bw^{\indk+1,i} = \bw^{\indk,i} + \bx^{\indk+1,i}(T) - \by^{\indk+1,i}$.
\end{enumerate}

The ADMM algorithm terminates when the number of iteration exceeds the maximal number of iteration or the following inequality holds
\begin{equation*}
\max\{\|\bx^{\indk+1,i}(T) - \bx^{\indk,i}(T)\|, \|\by^{\indk+1,i} - \by^{\indk,i}\|, \|\bw^{\indk+1,i} - \bw^{\indk,i}\|\}\leq \epsilon
\end{equation*}
for some positive threshold $\epsilon$. 
If the ADMM algorithm for the $i$-th subproblem terminates at the $N_i$-th step, we get the output parameters $\bar \bw^i$ and $\bar \by^i$ by
\begin{equation*}
\bar \bw^i = \bw^{N_i, i}, \quad  \bar \by^i = \by^{N_i, i}.
\end{equation*}
Then, the solution $\V$ to the HJ PDE with the general terminal condition in~\eqref{eqt:J_min_Ji} is computed using the neural network architecture in Fig.~\ref{fig: nn_Riccati}, where the coefficients in the $i$-th neuron, denoted by $\Sxx_i(t)$, $\Sx_i(t)$ and $\Sc_i(t)$, are the solutions to the FVPs~\eqref{eqt: odeP},~\eqref{eqt: odeq} and~\eqref{eqt: oder} with terminal condition
\begin{equation*}
    \Sxx_i(T) = \rho I_n, \quad \Sx_i(T) = \rho(\bar\bw^{i} - \bar\by^{i}), \quad \Sc_i(T) = \frac{\rho}{2}\|\bar\bw^{i}-\bar\by^{i}\|^2.
\end{equation*}

\section{Implementations of the abstract neural network architectures}\label{sec:implementation}

In the neural network architectures depicted in Figs.~\ref{fig: nn_Riccati} and~\ref{fig: nn_u}, each neuron involves the functions $P_i(t), \bq_i(t)$ and $r_i(t)$, which are the solutions to the FVPs~\eqref{eqt: odeP},~\eqref{eqt: odeq} and~\eqref{eqt: oder}. As a result, the neural network architectures require a numerical solver for solving the matrix Riccati FVP~\eqref{eqt: odeP} and the FVPs~\eqref{eqt: odeq} and~\eqref{eqt: oder}. 

In the literature, there are many numerical methods developed for solving the matrix Riccati differential equation. In order to solve this equation with a general initial or terminal condition, different fundamental solutions are proposed, including but not limited to Davison-Maki fundamental solution \cite{Davison1973numerical,Kenney1985numerical}, symplectic fundamental solution \cite{Levin1959matrix}, and min-plus fundamental solution \cite{Mceneaney2008new, Dower2015new, Dower2015max, Deshpande2011maxplus}. 
Recently, there are some non-traditional methods developed for solving Riccati equations using ant colony programming \cite{kamali2015study}, genetic programming \cite{Balasubramaniam2009solution} and neural networks \cite{Balasubramaniam2006solution,Samath2010solution}.

\begin{algorithm}[htbp]
\SetAlgoLined
\SetKwInOut{Input}{Input}
\SetKwInOut{Output}{Output}
\Input{Time $t\in[t_0,T)$ and the step number $N$.}
\Output{The solution $z(t)$ at time $t$.}
 Initialization: set $z_N \doteq z_T$, $t_N\doteq T$ and $\Delta t\doteq \frac{T-t}{N}$\;
 \For{$k = N, N-1, \dots, 1$}{
 $\delta_1 \doteq -\Delta t  g(t_k,z_k)$ and $w_1 \doteq z_k + \frac{1}{2}\delta_1$\;
 $\delta_2 \doteq -\Delta t  g(t_k-\frac{\Delta t}{2},w_1)$ and $w_2 \doteq z_k + \frac{1}{2}\delta_2$\;
 $\delta_3 \doteq -\Delta t  g(t_k - \frac{\Delta t}{2},w_2)$ and $w_3 \doteq z_k + \delta_3$\;
 $\delta_4 \doteq -\Delta t  g(t_k-\Delta t,w_3)$\;
 Update $z_{k-1}\doteq z_k + \frac{\delta_1}{6}+ \frac{\delta_2}{3}+ \frac{\delta_3}{3}+ \frac{\delta_4}{6}$ and 
 $t_{k-1} \doteq t_k - \Delta t$\;
 }
 The output $z(t)$ is given by $z_0$.
 \caption{The fourth order Runge-Kutta method for solving the FVP~\eqref{eqt: general_backward_ODE}. \label{alg:RK4}}
\end{algorithm}

\begin{figure}[htbp]
\centering
\includegraphics[width = \textwidth]{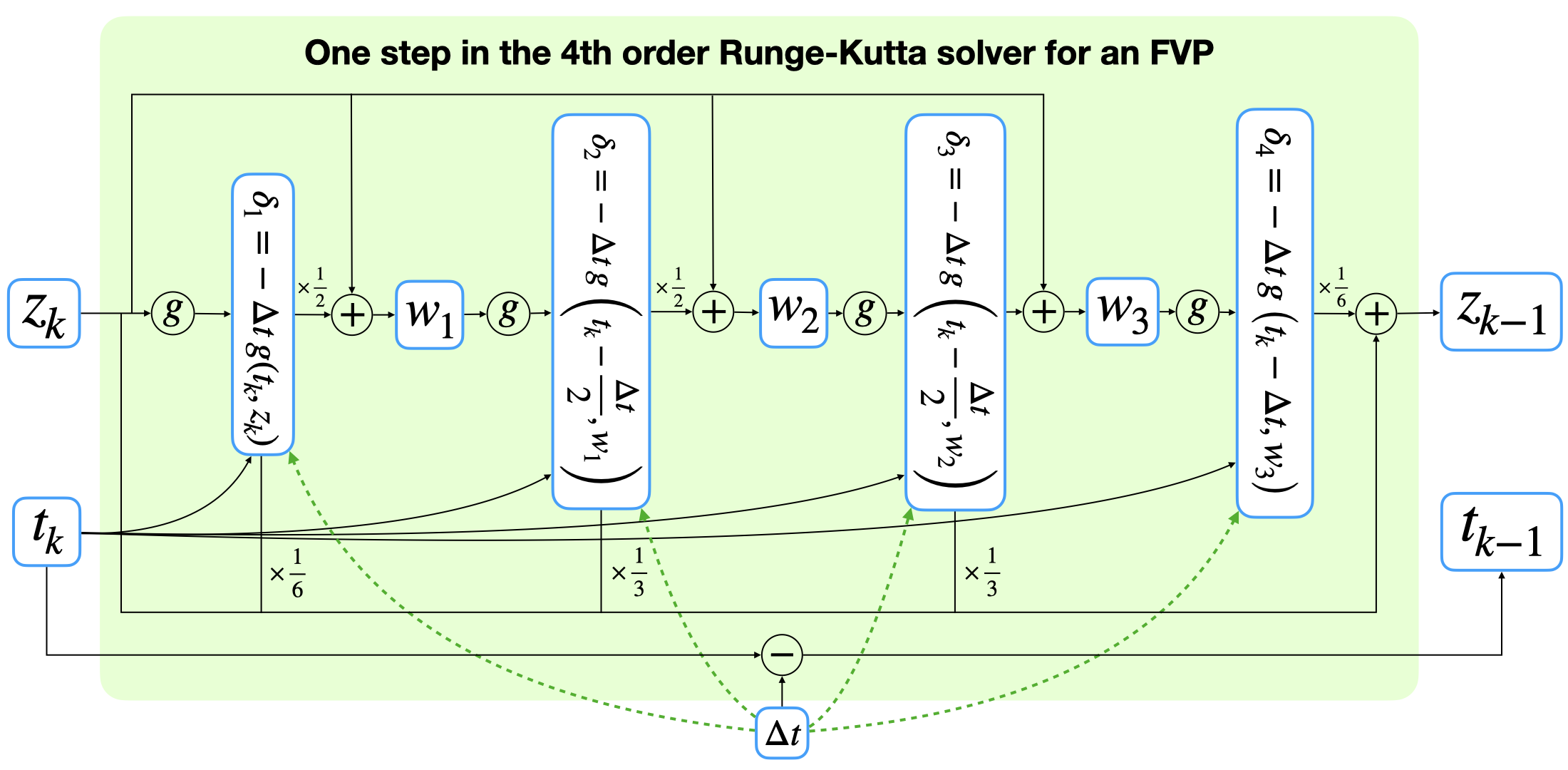}
\caption{
Illustration of the Resnet architecture that represents one step in the fourth order Runge-Kutta solver shown in Algorithm~\ref{alg:RK4} for solving a general FVP~\eqref{eqt: general_backward_ODE}.}
\label{fig: nn_RK4}
\end{figure}

\begin{figure}[htbp]
\centering
\includegraphics[width = \textwidth]{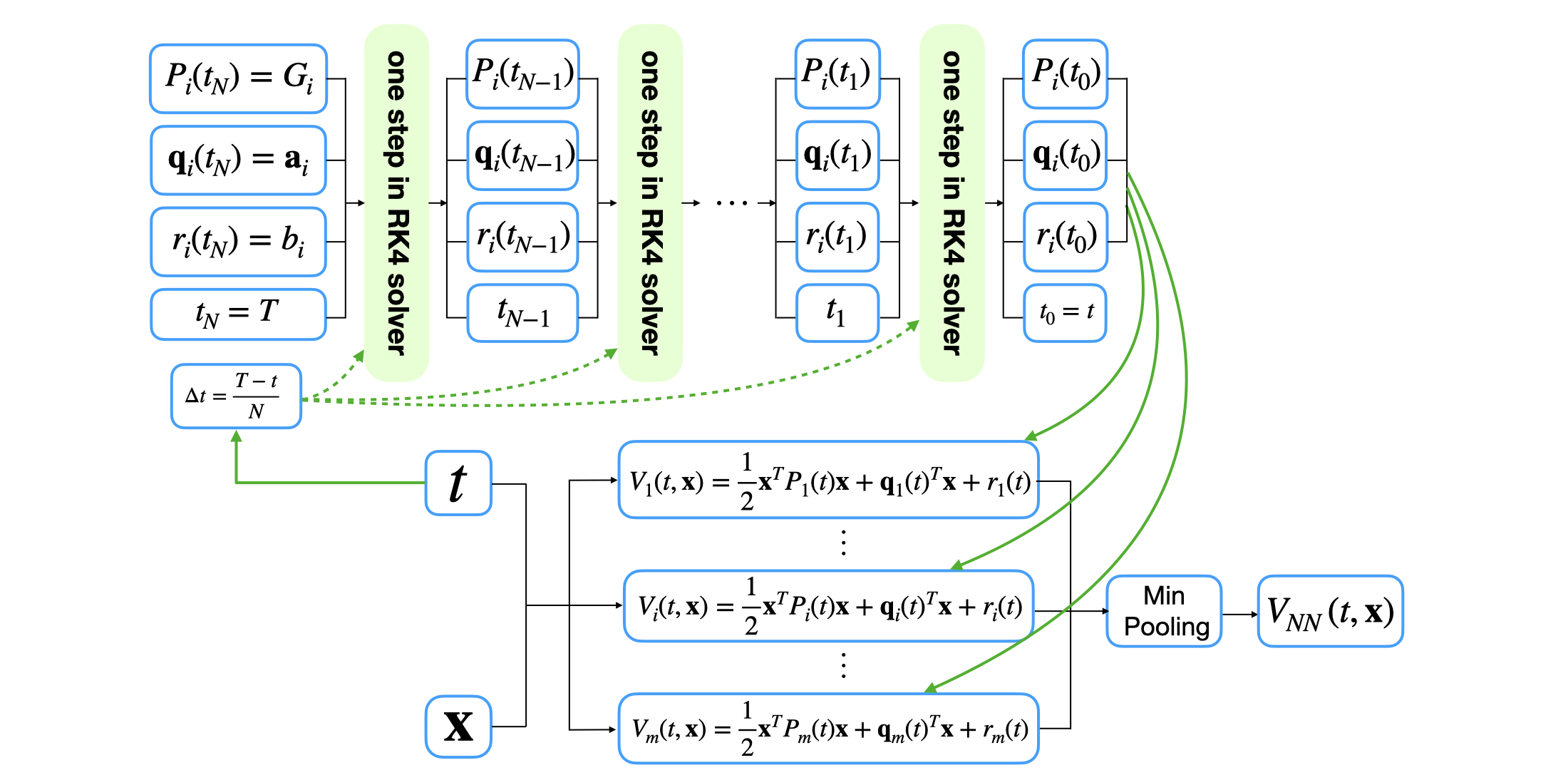}
\caption{
An implementation of the abstract neural network architecture defined by~\eqref{eqt: nn_u} where $\{V_i(t,\bx)\}_{i=1}^m$ are computed using the Resnet neural network depicted in Fig.~\ref{fig: nn_RK4}.
}
\label{fig: nn_V_RK4}
\end{figure}

\begin{figure}[htbp]
\centering
\includegraphics[width = \textwidth]{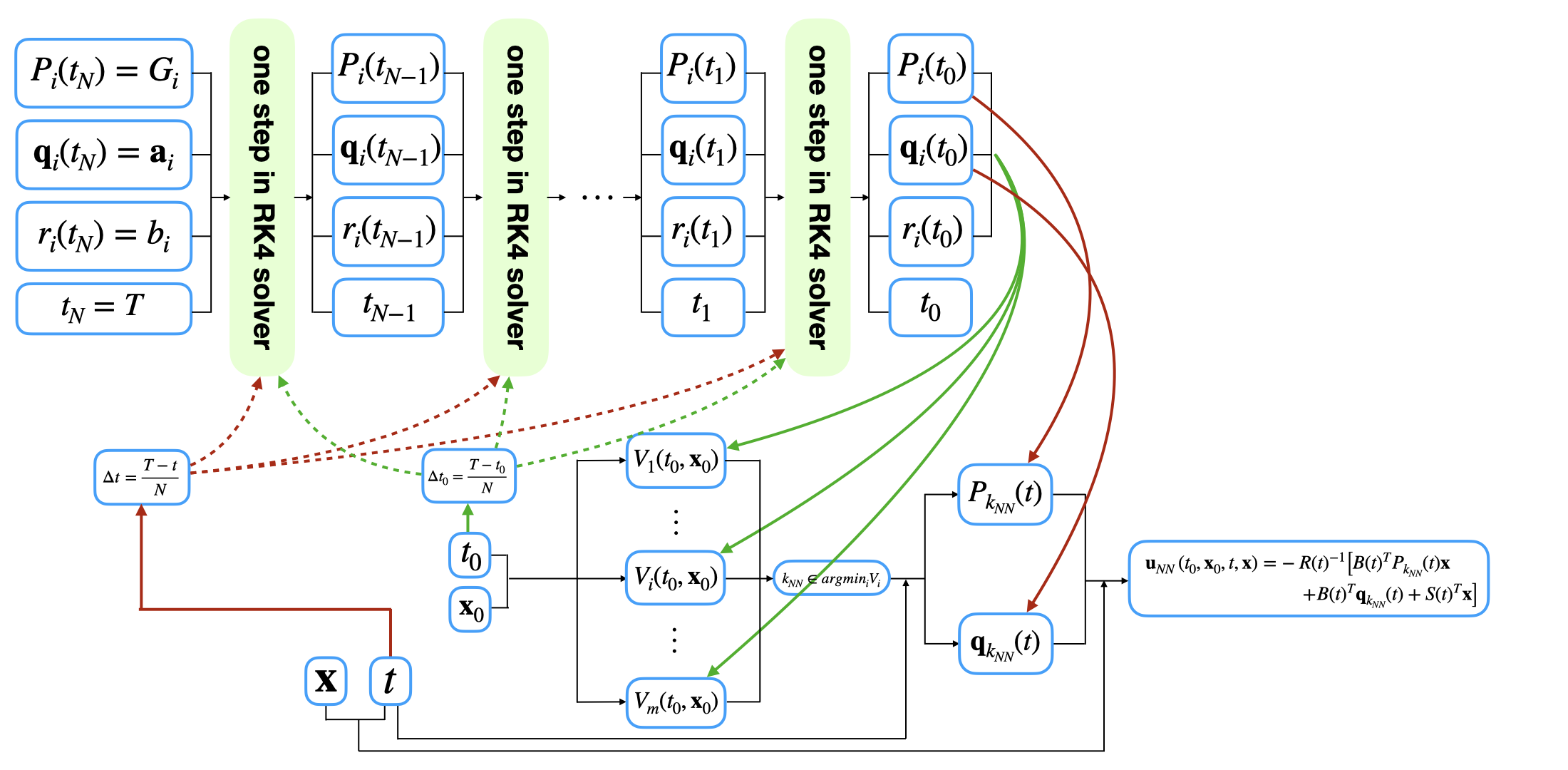}
\caption{
An implementation of the abstract neural network architecture defined by~\eqref{eqt: nn_u} where the neurons $\{V_i(\tzero, \xzero)\}_{i=1}^m$, $\bq_{\knn}(t)$ and $P_{\knn}(t)$ (where $\knn$ is the index defined by~\eqref{eqt: defk}) are computed using a fourth order Runge-Kutta method depicted in Fig.~\ref{fig: nn_RK4}.}
\label{fig: nn_u_RK4}
\end{figure}

We adopt the fourth order Runge-Kutta method to solve the general FVP
\begin{equation}\label{eqt: general_backward_ODE}
\begin{dcases}
\dot{z}(t) = g(t,z(t)) & t\in[t_0,T],\\
z(T) = z_T,
\end{dcases}
\end{equation}
where the function $z\colon [t_0,T]\to \R^\alpha$ (for a positive integer $\alpha$) is an absolutely continuous function solving the FVP almost everywhere, and the source term $g\colon [t_0,T]\times \R^\alpha\to\R^\alpha$ is continuous with respect to $t$ and uniformly Lipschitz with respect to $z$. The fourth order Runge-Kutta algorithm for solving this FVP is reviewed in Algorithm~\ref{alg:RK4}.
Note that the Runge-Kutta solver can be expressed using a neural network architecture (see, for instance,~\cite{Anastassi2014constructing}). For illustration, we show in Fig.~\ref{fig: nn_RK4} the architecture corresponding to one step of the fourth order Runge-Kutta solver. 
This architecture belongs to the class of Resnet architectures proposed in~\cite{He2016Deep}.
With this connection, the Runge-Kutta solver can be implemented using standard neural network languages, which can be converted to executable codes on the dedicated hardware designed for neural networks. By employing the Runge-Kutta solver to evaluate each abstract neuron in the proposed architectures, we obtain their implementations using Resnet-type deep neural networks.
The illustrations of the deep neural network implementations for $\VNN$ and $\unn$ are shown in Figs.~\ref{fig: nn_V_RK4} and~\ref{fig: nn_u_RK4}, respectively.
From the neural network function $\unn$, we compute the optimal trajectory $\bx^*$ by solving~\eqref{eqt: ODE_optx} using the fourth order Runge-Kutta method. Then, we get the open loop control $\bu^*$ as described in Remark~\ref{rem:compute_control}. The TensorFlow implementations of these architectures for our examples are given in \url{https://github.com/TingweiMeng/NN_HJ_minplus}.

\begin{remark}
It appears that the residual grows with $\|\bx\|^2$ for a fixed time $t\in[0,T)$. It is expected because it follows from the Runge-Kutta error estimation (see, for instance, \cite{butcher2016numerical}). For a fixed time step (i.e., a fixed number of layers in our proposed neural network in Fig.~\ref{fig: nn_V_RK4}), the error at $(\bx,t)$ for a fixed time variable $t\in[0,T)$ and any spatial variable $\bx\in \Omega$ is bounded in a compact set $\Omega$, but not in the whole domain. However, this error will converge to zero as the number of layers goes to infinity.
\end{remark}

Our implementations are based on the fourth order Runge-Kutta solver. Note that other ODE solvers can also be applied to compute the neurons in the abstract neural network architectures in Figs.~\ref{fig: nn_Riccati} and~\ref{fig: nn_u}. Each ODE solver which can be represented using neural network architectures provides possible neural network implementations for the two abstract architectures in Figs.~\ref{fig: nn_Riccati} and~\ref{fig: nn_u}, which therefore provide possibilities for leveraging different neural network architectures to solve high dimensional HJ PDEs~\eqref{eqt: HJ_quadratic_terminal} and corresponding optimal control problems~\eqref{eqt: optctrl_quadratic}.

We will show three numerical experiments. In these experiments, we assume the coefficients $\Cpp,\Cxx,\Cxp,\Muu,\Mxx,\Mxu, \Af,\Bf$ satisfy the assumptions (A1)-(A2). The first example is shown in Sec.~\ref{subsec: test_const}, which has constant coefficients. The second example is shown in Sec.~\ref{subsec: test_tdep}, whose coefficients depend on the time variable. And the third example is shown in Sec.~\ref{subsec: test_newton}, which is a slightly modified version of the HJ PDE~\eqref{eqt: HJ_quadratic_terminal} and the optimal control problem~\eqref{eqt: optctrl_quadratic} considered in this paper.

In each example, we use the deep Resnet implementation depicted in Fig.~\ref{fig: nn_V_RK4}
to solve the viscosity solution to the HJ PDEs and the value function in the optimal control problems at different time $t$.
For different terminal time $T$, we solve the corresponding optimal controls and optimal trajectories with different initial position $\bx_0$ by the method described in Remark~\ref{rem:compute_control}
using the deep Resnet neural network implementation depicted in Fig.~\ref{fig: nn_u_RK4}.
If not mentioned explicitly, we use $40$ Runge-Kutta layers to compute the viscosity solution $\VNN$ and $400$ Runge-Kutta layers to compute the optimal controls and optimal trajectories.

To show the solution $\VNN$ in high dimensional cases, we plot two dimensional slices of the function $\bx\mapsto \VNN(t, \bx)$ for different time $t$.
We consider the points $\bx = (x_1, x_2, \bzero)\in\Rn$ where $(x_1,x_2)\in\R^2$ is any grid point in a two dimensional rectangular domain and $\bzero$ denotes the zero vector in $\R^{n-2}$. In each figure, the color is given by the function value $\VNN(t,x_1,x_2, \bzero)$, and the x and y axes represent the variables $x_1$ and $x_2$, respectively. 
To show the errors of the viscosity solution, we compute the maximal absolute value of the residual (i.e., $\max_{i\in\{1,\dots,m\}} |-\frac{\partial \V_i}{\partial t} + H(t,\bx,\nabla_{\bx}\V_i)|$) in each example, where $\V_i$ is defined in~\eqref{eqt: nn_S}. 
We set the number of Runge-Kutta layers to be $20$, $40$, and $80$ to show the dependence of error on the number of layers. The residual values in different example, at different times $t$, and computed using different number of Runge-Kutta layers are shown in Table~\ref{tab:err_table}. From the table, we observe that the magnitude of the absolute values of the residuals is in general small (less than $10^{-6}$), which provides a numerical validation that each $\V_i$ approximately satisfies the differential equation in~\eqref{eqt: HJ_quadratic_terminal}. Since the solution operator of the HJ PDE~\eqref{eqt: HJ_quadratic_terminal} is linear with respect to the min-plus algebra, our proposed deep neural network architecture in Figure~\ref{fig: nn_V_RK4} indeed approximates the viscosity solution to the HJ PDE. The errors also decrease as the number of layers goes to infinity. This observation validates the error analysis of Runge-Kutta solvers (see, for instance, \cite{butcher2016numerical}).

\begin{table}[htbp]
    \centering
    \begin{tabular}{c|c|c|c|c}
        \hline
         & \# RK layers & $t=0.25$ & $t=0.5$ & $t=0.75$ \\
        \hline
        \hline
        \multirow{3}{*}{\shortstack{Example 1\\ in Section~\ref{subsec: test_const}}
        }
         & 20 &6.97E-06  & 1.45E-06 &6.59E-08 \\
        & 40 & 4.21E-07 & 8.93E-08 & 4.12E-09\\
        & 80 & 2.59E-08 & 5.54E-09 & 2.57E-10\\
        \hline
        \multirow{3}{*}{\shortstack{Example 2\\ in Section~\ref{subsec: test_tdep}}}
         & 20 & 1.29E-07 & 3.28E-08 & 2.91E-09\\
        & 40 & 7.94E-09 & 2.02E-09 & 1.80E-10\\
        & 80 & 4.92E-10 & 1.25E-10 & 1.12E-11\\
        \hline
        \multirow{3}{*}{\shortstack{Example 4\\ in Section~\ref{subsec:example_admm}}}
         & 20 & 2.24E-07 & 4.49E-08 & 3.08E-09\\
        & 40 & 1.39E-08 & 2.77E-09 & 1.90E-10\\
        & 80 & 8.64E-10 & 1.72E-10 & 1.18E-11\\
        \hline
        \hline
        & \# RK layers & $t=0.75$ & $t=0.95$ & $t=0.995$ \\
        \hline \hline
        \multirow{3}{*}{\shortstack{Example 3\\ in Section~\ref{subsec: test_newton}}}
         & 20 & 3.95E-04 & 1.21E-04 & 1.44E-08\\
        & 40 & 2.14E-05 & 7.10E-06 & 9.01E-10\\
        & 80 & 1.25E-06 & 4.30E-07 & 5.67E-11\\
        \hline
    \end{tabular}
    \hfill 
    \caption{We show the maximal absolute residual $\max_{i\in\{1,\dots,m\}} |-\frac{\partial \V_i}{\partial t} + H(t,\bx,\nabla_{\bx}\V_i)|$ in each example, where $\V_i$ is defined in~\eqref{eqt: nn_S}. The residual values are computed at different time~$t$, with different Runge-Kutta (RK) layers (which is related to the number of layers in the proposed neural networks).}
    \label{tab:err_table}
\end{table}

Also, to illustrate the optimal controls and optimal trajectories in each example, we consider different initial positions $\bx_0=(x,\bzero)\in\Rn$ (where $x$'s are the grid points in a one-dimensional interval and $\bzero$ denotes the zero vector in $\R^{n-1}$) and a fixed initial time $\tzero=0$. To avoid ambiguity, we use the notations $\bu^*(s|\bx_0)$ and $\bx^*(s|\bx_0)$ to denote the optimal control and trajectory at time $s$ with initial position $\bx_0$. For each initial time $\tzero=0$ and initial position $\bx_0\in\Rn$, we compute the corresponding optimal control, denoted by $[0,T]\ni s\mapsto \bu^*(s|\bx_0)\in \R^l$, and the corresponding optimal trajectory, denoted by $[0,T]\ni s\mapsto \bx^*(s|\bx_0)\in\Rn$. For the one dimensional problem, we choose each initial position $x_0\in\R$ to be a grid point in a one-dimensional interval, and plot the graphs of $u^*(\cdot|x_0)$ or $x^*(\cdot|x_0)$ with different $x_0$ in one figure. For the high dimensional problem, we choose different initial positions $\bx_0 = (x, \bzero)$ where $x$'s are the grid points in a one-dimensional interval and $\bzero$ denotes the zero vector in $\R^{n-1}$, and solve the optimal controls $\bu^*(\cdot|\bx_0) = (u_1^*(\cdot|\bx_0), \dots, u_l^*(\cdot|\bx_0))$ and the optimal trajectories $\bx^*(\cdot|\bx_0) = (x_1^*(\cdot|\bx_0),\cdots, x_n^*(\cdot|\bx_0))$. Then, the graphs of the some components of the optimal controls $\bu^*(\cdot|\bx_0)$ or optimal trajectories $\bx^*(\cdot|\bx_0)$ are plotted in each figure.

\subsection{An optimal control problem with constant coefficients} \label{subsec: test_const}

\begin{figure}[htbp]
\begin{minipage}[b]{.49\linewidth}
  \centering
  \includegraphics[width=.99\textwidth]{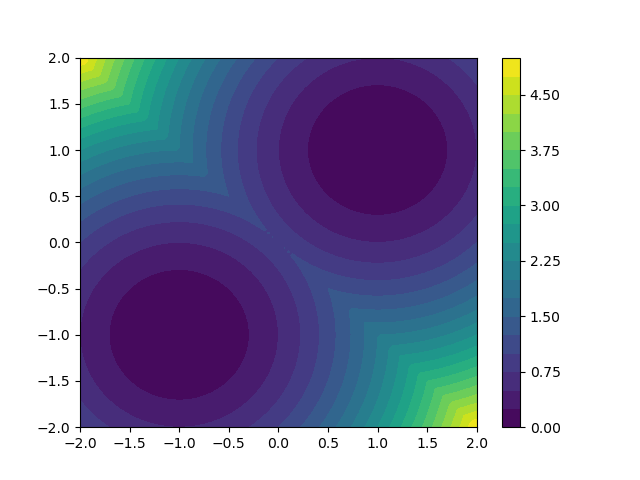}
  \centerline{\footnotesize{(a) $t=1$}}\medskip
\end{minipage}
\hfill  
\begin{minipage}[b]{.49\linewidth}
  \centering
  \includegraphics[width=.99\textwidth]{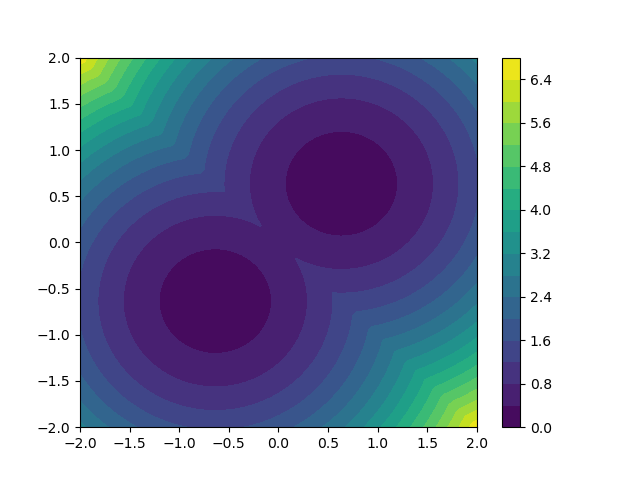}
  \centerline{\footnotesize{(b) $t=0.75$}}\medskip
\end{minipage}
\begin{minipage}[b]{.49\linewidth}
  \centering
    \includegraphics[width=.99\textwidth]{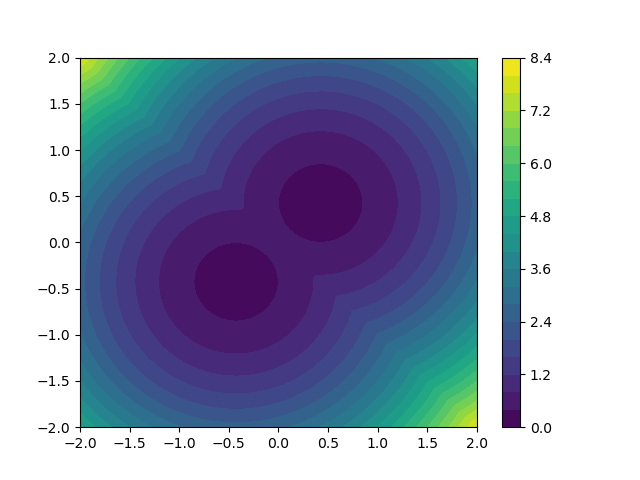}
  \centerline{\footnotesize{(c) $t=0.5$}}\medskip
\end{minipage}
\hfill
\begin{minipage}[b]{.49\linewidth}
  \centering
    \includegraphics[width=.99\textwidth]{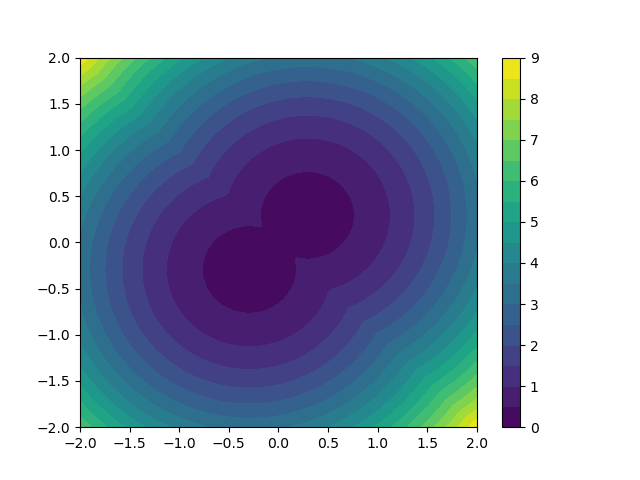}
  \centerline{\footnotesize{(d) $t=0.25$}}\medskip
\end{minipage}
\caption{The viscosity solution $\VNN$ to the $16$ dimensional HJ PDE~\eqref{eqt: HJ_quadratic_terminal} with Hamiltonian~\eqref{eqt:eg1_H}, terminal data~\eqref{eqt: J_testconst_16d} and terminal time $T=1$ is computed using the proposed abstract neural network architecture~\eqref{eqt: nn_S} (depicted in Fig.~\ref{fig: nn_Riccati}) with the implementation depicted in Fig.~\ref{fig: nn_V_RK4}. 
The two dimensional slices of $\VNN$ at time $t=1$ (i.e., the terminal cost), $t=0.75$, $t=0.5$ and $t=0.25$ are shown in the subfigures (a), (b), (c) and (d), respectively.
The color in each subfigure shows the solution value $\VNN(t, \bx)$, where the spatial variable $\bx$ is in the form of $(x_1,x_2,\bzero)\in\R^{16}$ (with $\bzero$ denoting the zero vector in $\R^{14}$) for some points $x_1\in\R$ and $x_2\in\R$ which are represented by the $x$ and $y$ axes.
\label{fig: test_const_16d}}
\end{figure}

\begin{figure}[htbp]
\begin{minipage}[b]{.49\linewidth}
  \centering
  \includegraphics[width=.99\textwidth]{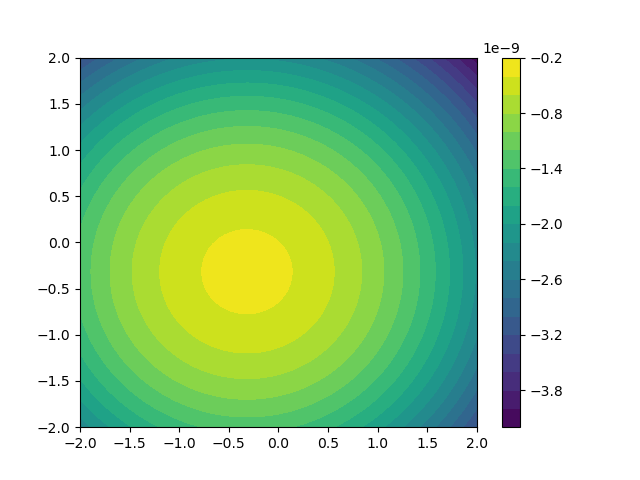}
  \centerline{\footnotesize{(a) $t=0.75$}}\medskip
\end{minipage}
\hfill  
\begin{minipage}[b]{.49\linewidth}
  \centering
  \includegraphics[width=.99\textwidth]{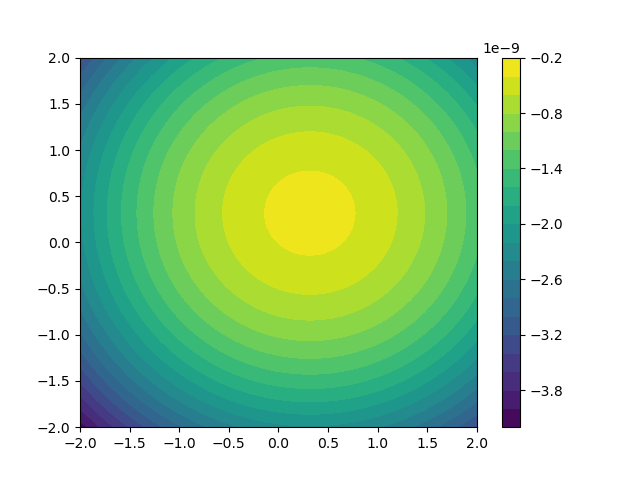}
  \centerline{\footnotesize{(b) $t=0.75$}}\medskip
\end{minipage}
\begin{minipage}[b]{.49\linewidth}
  \centering
    \includegraphics[width=.99\textwidth]{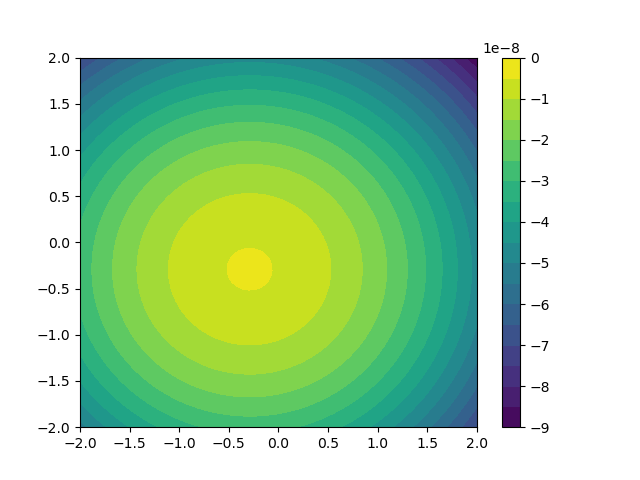}
  \centerline{\footnotesize{(c) $t=0.5$}}\medskip
\end{minipage}
\hfill
\begin{minipage}[b]{.49\linewidth}
  \centering
    \includegraphics[width=.99\textwidth]{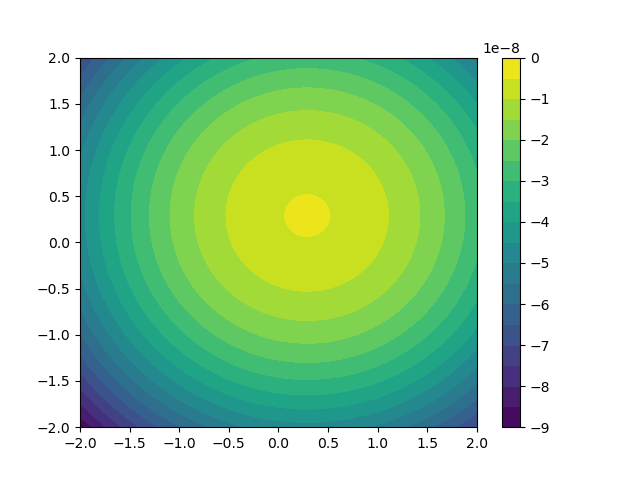}
  \centerline{\footnotesize{(d) $t=0.5$}}\medskip
\end{minipage}
\begin{minipage}[b]{.49\linewidth}
  \centering
    \includegraphics[width=.99\textwidth]{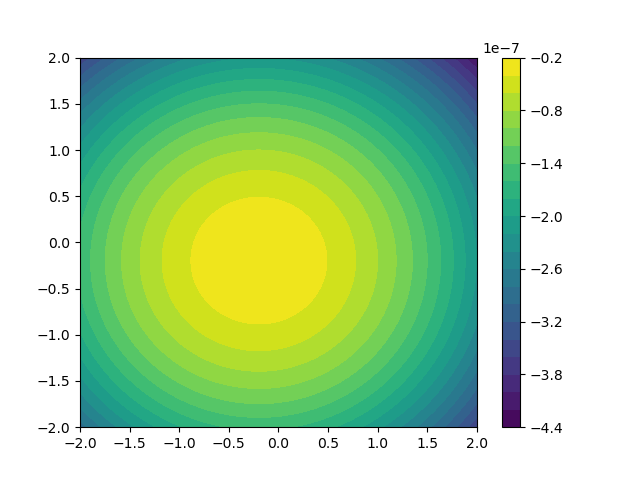}
  \centerline{\footnotesize{(e) $t=0.25$}}\medskip
\end{minipage}
\hfill
\begin{minipage}[b]{.49\linewidth}
  \centering
    \includegraphics[width=.99\textwidth]{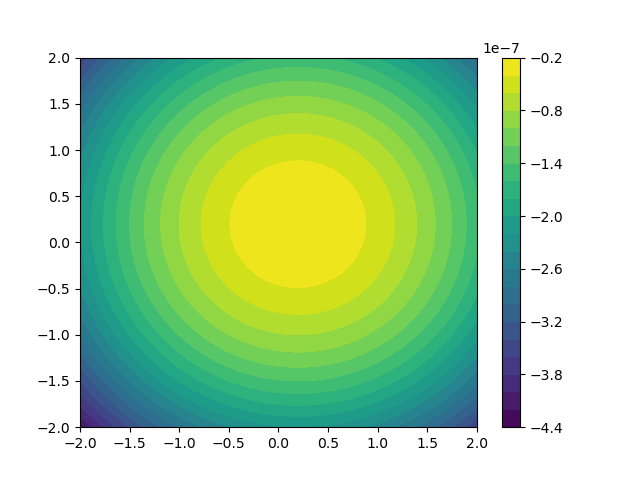}
  \centerline{\footnotesize{(f) $t=0.25$}}\medskip
\end{minipage}
\caption{The residual $-\frac{\partial \V_i}{\partial t} + H(t,\bx,\nabla_{\bx}\V_i)$ in the HJ PDE~\eqref{eqt: HJ_quadratic_terminal} with Hamiltonian~\eqref{eqt:eg1_H} (with terminal time $T=1$) is shown for each $i\in\{1,\dots,m\}$ at different time $t$, where $\V_i$ is defined in~\eqref{eqt: nn_S}. Figures (a), (c), (e) show the residuals for $\V_1$ at time $t=0.75$, $t=0.5$ and $t=0.25$, while figures (b), (d), (f) show the residuals for $\V_2$ at time $t=0.75$, $t=0.5$ and $t=0.25$, respectively. 
In each subfigure, we show the two dimensional slices of the residual function. The color shows the residual value at $(t, \bx)$, where the spatial variable $\bx$ is in the form of $(x_1,x_2,\bzero)\in\R^{16}$ (with $\bzero$ denoting the zero vector in $\R^{14}$) for some points $x_1\in\R$ and $x_2\in\R$ which are represented by the $x$ and $y$ axes.
\label{fig:eg1_pde_err}}
\end{figure}


\begin{figure}[htbp]
\begin{minipage}[b]{.49\linewidth}
  \centering
  \includegraphics[width=.99\textwidth]{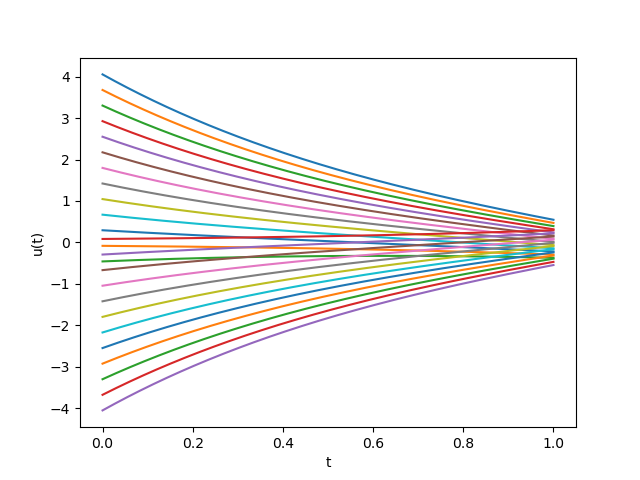}
  \centerline{\footnotesize{(a) controls, $T=1$}}\medskip
\end{minipage}
\hfill  
\begin{minipage}[b]{.49\linewidth}
  \centering
  \includegraphics[width=.99\textwidth]{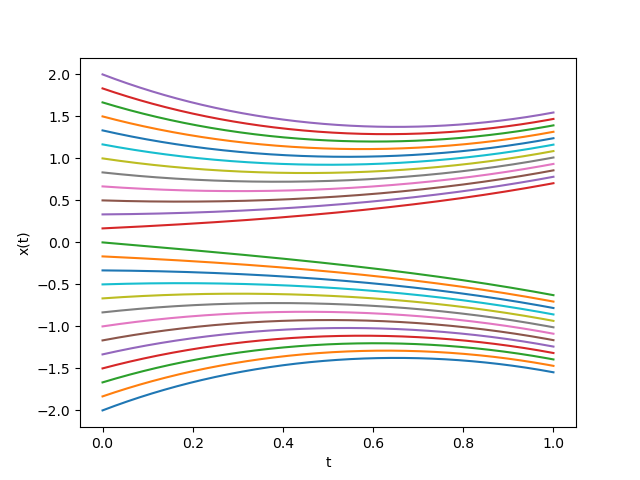}
  \centerline{\footnotesize{(b) trajectories, $T=1$}}\medskip
\end{minipage}
\begin{minipage}[b]{.49\linewidth}
  \centering
  \includegraphics[width=.99\textwidth]{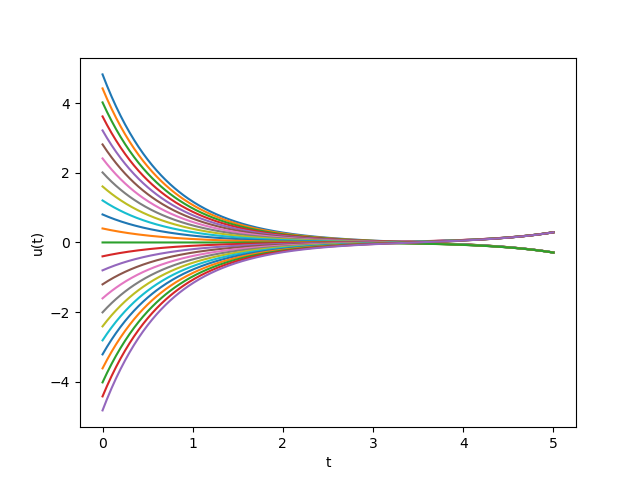}
  \centerline{\footnotesize{(c) controls, $T=5$}}\medskip
\end{minipage}
\hfill
\begin{minipage}[b]{.49\linewidth}
  \centering
    \includegraphics[width=.99\textwidth]{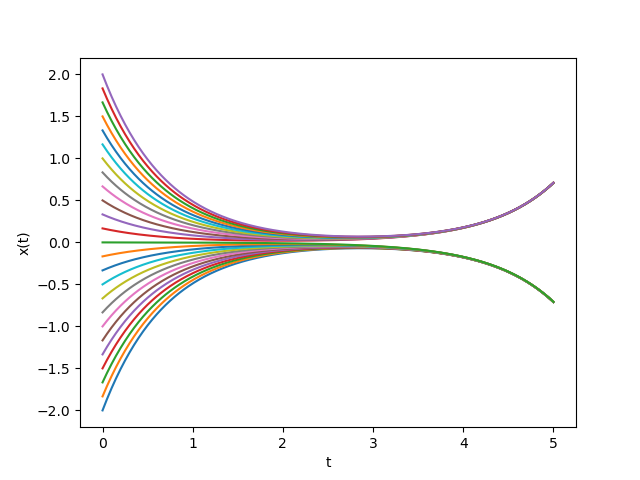}
  \centerline{\footnotesize{(d) trajectories, $T=5$}}\medskip
\end{minipage}
\caption{The open loop optimal controls and the corresponding optimal trajectories in the $16$ dimensional optimal control problem~\eqref{eqt: optctrl} with Lagrangian~\eqref{eqt:eg1_L}, source term~\eqref{eqt:eg1_f}, terminal cost~\eqref{eqt: J_testconst_16d} and different terminal time $T=1$, $5$ are computed using the proposed abstract neural network architecture~\eqref{eqt: nn_u} with the implementation depicted in Fig.~\ref{fig: nn_u_RK4}. Several graphs of the first component of the optimal controls with $T=1$ are shown in (a), and the first component of the corresponding optimal trajectories are shown in (b). Several graphs of the first component of the optimal controls with $T=5$ are shown in (c), and the first component of the corresponding optimal trajectories are shown in (d).
\label{fig: test_const_16dxu}}
\end{figure}

We consider the optimal control problem~\eqref{eqt: optctrl_quadratic} with the following constant coefficients
\begin{equation*}
\begin{dcases}
    l = n,\\
    \Muu = \Mxx = \Cpp = \Cxx = \Cxp = \Af = \Bf = I_n,\\
    \Mxu = \On,
\end{dcases}
\end{equation*}
where $I_n$ denotes the identity matrix in $\matnn$ and $\On$ denotes the zero matrix in $\matnn$. 
With these coefficients, we solve the optimal control problem~\eqref{eqt: optctrl} whose Lagrangian $L$ in~\eqref{eqt:def_quad_L_f} is defined by
\begin{equation}\label{eqt:eg1_L}
    L(t,\bx,\bu) =  \frac{1}{2}\|\bx\|^2 +  \frac{1}{2}\|\bu\|^2  \quad \forall t\in[0,T], \bx\in\Rn, \bu\in\R^l,
\end{equation}
and the source term $f$ in~\eqref{eqt:def_quad_L_f} is defined by
\begin{equation}\label{eqt:eg1_f}
    f(t,\bx,\bu) = \bx+\bu \quad \forall t\in[0,T], \bx\in\Rn, \bu\in\R^l.
\end{equation}
The corresponding HJ PDE is in the form of~\eqref{eqt: HJ_quadratic_terminal} where the Hamiltonian $H$ is defined by
\begin{equation}\label{eqt:eg1_H}
    H(t,\bx,\bp)= \frac{1}{2}\| \bp\|^2 - \frac{1}{2}\|\bx\|^2 - \langle \bp, \bx\rangle
    \quad \forall t\in[0,T], \bx,\bp\in\Rn.
\end{equation}
With these coefficients, the differential equations for $P_i$, $\bq_i$ and $r_i$ read
\begin{equation*}
    \begin{split}
    \dot{\Sxx}_i(t) &=  \Sxx_i(t)^T\Sxx_i(t) - 2\Sxx_i(t) - I_n, \\
    \dot{\Sx}_i(t) &= \Sxx_i(t)^T\Sx_i(t) - \Sx_i(t),\\
    \dot{\Sc}_i(t) &= \frac{1}{2}\|\Sx_i(t)\|^2,
    \end{split}
\end{equation*}
for each $t\in (0,T)$.
We consider this high dimensional problem with $n=16$ and $m=2$, where the terminal data $\J$ is defined by
    \begin{equation}\label{eqt: J_testconst_16d}
        \J(\bx) = \min\left\{ \frac{1}{2}\left(\sum_{i=1}^2(x_i +1)^2 + \sum_{i=3}^{16} x_i^2\right), \frac{1}{2}\left(\sum_{i=1}^2(x_i -1)^2 + \sum_{i=3}^{16} x_i^2\right)\right\},
    \end{equation}
    for each $\bx=(x_1,\dots,x_{16})\in \R^{16}$.

The viscosity solution to the HJ PDE~\eqref{eqt: HJ_quadratic_terminal} with Hamiltonian in~\eqref{eqt:eg1_H} and terminal data in~\eqref{eqt: J_testconst_16d}
is computed using the proposed abstract neural network~\eqref{eqt: nn_S} (depicted in Fig.~\ref{fig: nn_Riccati}) with the implementation depicted in Fig.~\ref{fig: nn_V_RK4}.
The two dimensional slices of the solution $\VNN$ to the $16$ dimensional problem with terminal time $T=1$ is shown in Fig.~\ref{fig: test_const_16d}. In this figure, the solution $\VNN(t,x_1,x_2,\bzero)$ (where $\bzero$ denotes the zero vector in $\R^{n-2}$) at time $t = 1$, $t = 0.75$, $t=0.5$ and $t = 0.25$ is shown in the subfigures (a), (b), (c) and (d), respectively.
Recall that the $x$ and $y$ axes represent the first component $x_1$ and second component $x_2$ of the spatial variable, respectively.
We also show in Figure~\ref{fig:eg1_pde_err} the (two-dimensional slices of) the residual $-\frac{\partial \V_i}{\partial t} + H(t,\bx,\nabla_{\bx}\V_i)$ for each $i\in\{1,\dots,m\}$, where $\V_i$ is defined in~\eqref{eqt: nn_S}. Note that the magnitude of the residuals is in general small (less than $10^{-6}$), which provides a numerical validation that each $\V_i$ approximately satisfies the differential equation in~\eqref{eqt: HJ_quadratic_terminal}. Since the solution operator of the HJ PDE~\eqref{eqt: HJ_quadratic_terminal} is linear with respect to the min-plus algebra, our proposed deep neural network architecture in Figure~\ref{fig: nn_V_RK4} indeed approximates the viscosity solution to the HJ PDE. It appears that the residual grows with $|\bx|^2$ for a fixed time $t\in[0,T)$. It is expected because it follows from the Runge-Kutta error estimation. For a fixed time step (i.e., a fixed number of layers in our proposed neural network in Fig.~\ref{fig: nn_V_RK4}), the error is bounded in a compact set, but not in the whole domain. However, this error will converge to zero as the number of layers goes to infinity.

We compute optimal controls and optimal trajectories with different initial positions $\xzero=(x,\bzero)\in\Rn$ (where $\bzero$ denotes the zero vector in $\R^{n-1}$) and a fixed initial time $\tzero = 0$. The first component of the optimal controls and trajectories is illustrated in Fig.~\ref{fig: test_const_16dxu}.
The optimal controls with terminal time $T = 1$ are shown in (a), the optimal trajectories with terminal time $T = 1$ are shown in (b), the optimal controls with terminal time $T = 5$ are shown in (c), and the optimal trajectories with terminal time $T = 5$ are shown in (d).
There does appear to be a turnpike phenomenon for the longer time horizon (see, for instance,~\cite{GRUNE2020Exponential,Zaslavski2015Turnpike} and the references in there).

\subsection{An optimal control problem with time dependent coefficients} \label{subsec: test_tdep}
We consider the optimal control problem~\eqref{eqt: optctrl_quadratic} whose coefficients depend on time. The coefficients are chosen to be
\begin{equation*}
    \begin{dcases}
        l = n,\\
        \Muu(t) = 4e^{-t}I_n,\, \Mxx(t) = \frac{e^{-t}}{2}I_n,\, \Mxu(t) \equiv e^{-t}I_n,\, \Af(t) \equiv \frac{1}{2}I_n,\, \Bf(t) \equiv I_n,\\
        \Cpp(t) = \frac{e^t}{4}I_n,\, \Cxp(t) = \frac{1}{4}I_n,\, \Cxx(t) = \frac{e^{-t}}{4}I_n,
    \end{dcases}
\end{equation*}
for each $t\in[0,T]$, 
where $I_n$ denotes the identity matrix in $\matnn$.
With these coefficients, we solve the optimal control problem~\eqref{eqt: optctrl} whose Lagrangian $L$ in~\eqref{eqt:def_quad_L_f} reads
\begin{equation}\label{eqt:eg2_L}
    L(t,\bx,\bu) =  \frac{e^{-t}}{4}\|\bx\|^2 +  2 e^{-t}\|\bu\|^2 + e^{-t}\bx^T\bu \quad \forall t\in[0,T], \bx, \bu\in\Rn,
\end{equation}
and the source term $f$ in~\eqref{eqt:def_quad_L_f} reads
\begin{equation}\label{eqt:eg2_f}
    f(t,\bx,\bu) = \frac{\bx}{2} + \bu \quad \forall t\in[0,T], \bx, \bu\in\Rn.
\end{equation}
The corresponding HJ PDE is in the form of~\eqref{eqt: HJ_quadratic_terminal} where the Hamiltonian is defined by
\begin{equation}\label{eqt:eg2_H}
    H(t,\bx,\bp) = 
    \frac{e^t}{8}\|\bp\|^2 - \frac{e^{-t}}{8}\|\bx\|^2 - \frac{1}{4}\bp^T \bx
    \quad \forall t\in[0,T], \bx,\bp\in\Rn.
\end{equation}
With these coefficients, the differential equations for $P_i$, $\bq_i$ and $r_i$ read
\begin{equation*}
\begin{split}
    \dot{\Sxx}_i(t) &= \frac{e^t}{4}\Sxx_i(t)^T\Sxx_i(t) - \frac{1}{2}\Sxx_i(t) - \frac{e^{-t}}{4}I_n, \\
    \dot{\Sx}_i(t) &= \frac{e^t}{4}\Sxx_i(t)^T\Sx_i(t) - \frac{1}{4}\Sx_i(t), \\
    \dot{\Sc}_i(t) &= \frac{e^{t}}{8}\|\Sx_i(t)\|^2,
\end{split}
\end{equation*}
for each $t\in (0,T)$.
The running cost~\eqref{eqt:eg2_L} involves a discount factor of $1$. If the terminal cost was similarly discounted, we
would expect to see this appear in the HJ PDE as a $-V$ term, see, e.g.,~\cite[Section~III.3.1]{Bardi1997Optimal}.

\begin{figure}[htbp]
\begin{minipage}[b]{.49\linewidth}
  \centering
  \includegraphics[width=.99\textwidth]{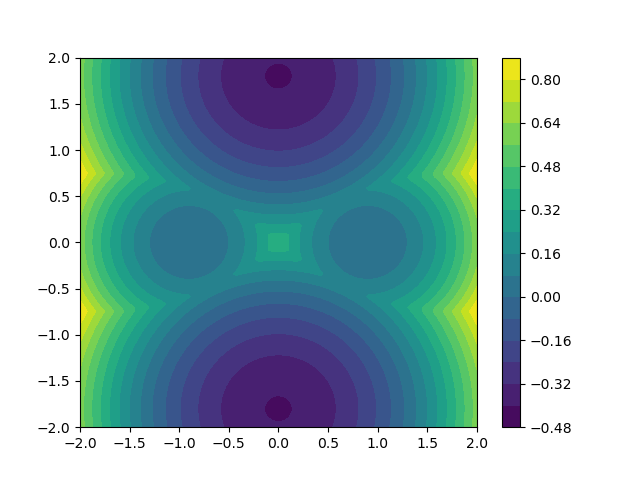}
  \centerline{\footnotesize{(a) $t=1$}}\medskip
\end{minipage}
\hfill  
\begin{minipage}[b]{.49\linewidth}
  \centering
  \includegraphics[width=.99\textwidth]{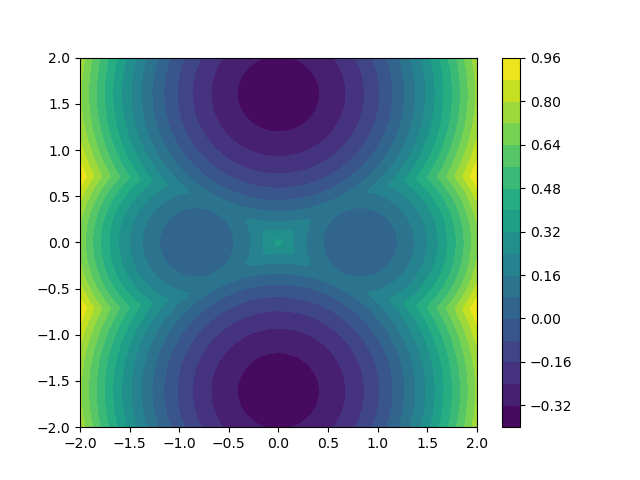}
  \centerline{\footnotesize{(b) $t=0.75$}}\medskip
\end{minipage}
\begin{minipage}[b]{.49\linewidth}
  \centering
    \includegraphics[width=.99\textwidth]{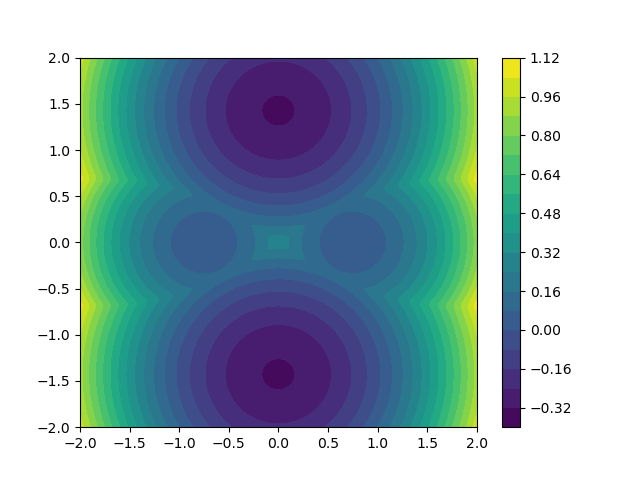}
  \centerline{\footnotesize{(c) $t=0.5$}}\medskip
\end{minipage}
\hfill
\begin{minipage}[b]{.49\linewidth}
  \centering
    \includegraphics[width=.99\textwidth]{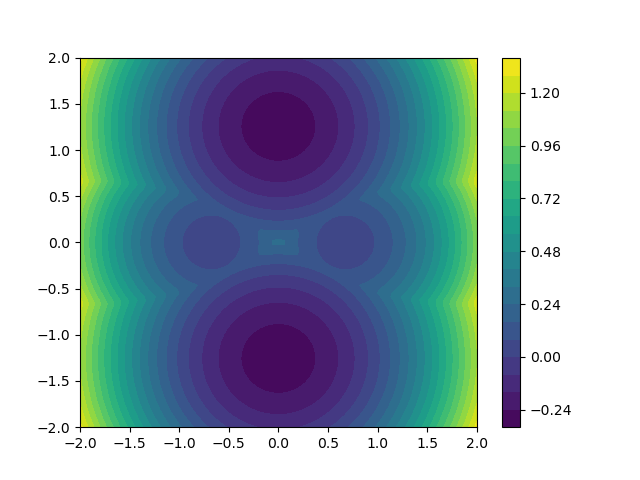}
  \centerline{\footnotesize{(d) $t=0.25$}}\medskip
\end{minipage}
\caption{The viscosity solution $\VNN$ to the $16$ dimensional HJ PDE~\eqref{eqt: HJ_quadratic_terminal} with Hamiltonian~\eqref{eqt:eg2_H}, terminal data~\eqref{eqt:J_min_Ji} (where $\J_i$'s are defined in~\eqref{eqt:Ji_tdep_16d}) and terminal time $T=1$ is computed using the proposed abstract neural network architecture~\eqref{eqt: nn_S} with the implementation depicted in Fig.~\ref{fig: nn_V_RK4}. The two dimensional slices of $\VNN$ at time $t=1$ (i.e., the terminal cost), $t=0.75$, $t=0.5$ and $t=0.25$ 
are shown in the subfigures (a), (b), (c) and (d), respectively.
The color in each subfigure shows the solution value $\VNN(t, \bx)$, where the spatial variable $\bx$ is in the form of $(x_1,x_2,\bzero)\in\R^{16}$ (where $\bzero$ is the zero vector in $\R^{14}$) for some points $x_1\in\R$ and $x_2\in\R$ which are represented by $x$ and $y$ axes.
\label{fig: test_tdep_16d}}
\end{figure}


\begin{figure}[htbp]
\begin{minipage}[b]{.49\linewidth}
  \centering
  \includegraphics[width=.99\textwidth]{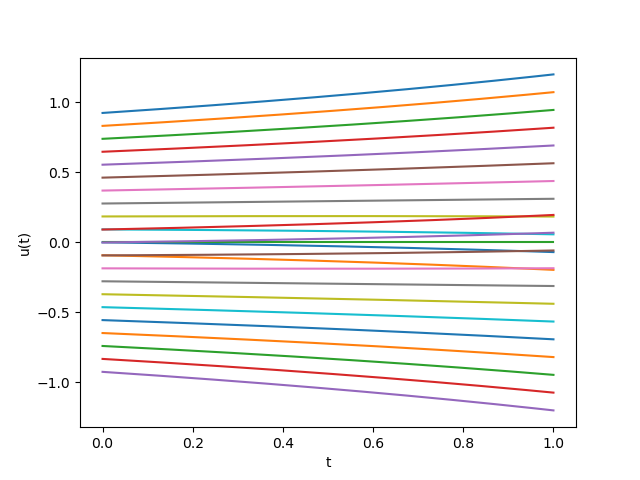}
  \centerline{\footnotesize{(a) controls, $T=1$}}\medskip
\end{minipage}
\hfill  
\begin{minipage}[b]{.49\linewidth}
  \centering
  \includegraphics[width=.99\textwidth]{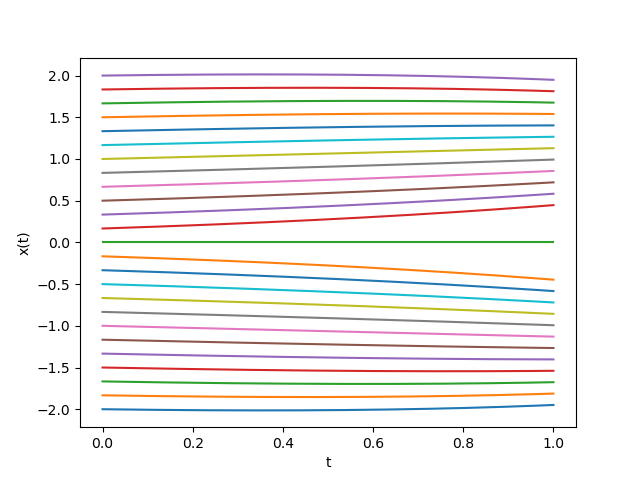}
  \centerline{\footnotesize{(b) trajectories, $T=1$}}\medskip
\end{minipage}
\begin{minipage}[b]{.49\linewidth}
  \centering
    \includegraphics[width=.99\textwidth]{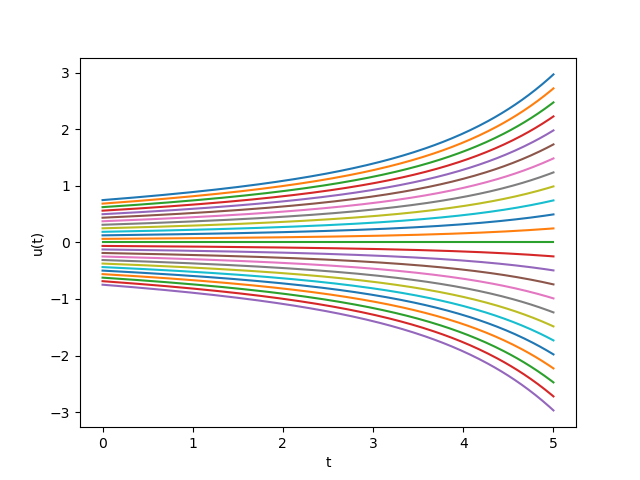}
  \centerline{\footnotesize{(c) controls, $T=5$}}\medskip
\end{minipage}
\hfill
\begin{minipage}[b]{.49\linewidth}
  \centering
    \includegraphics[width=.99\textwidth]{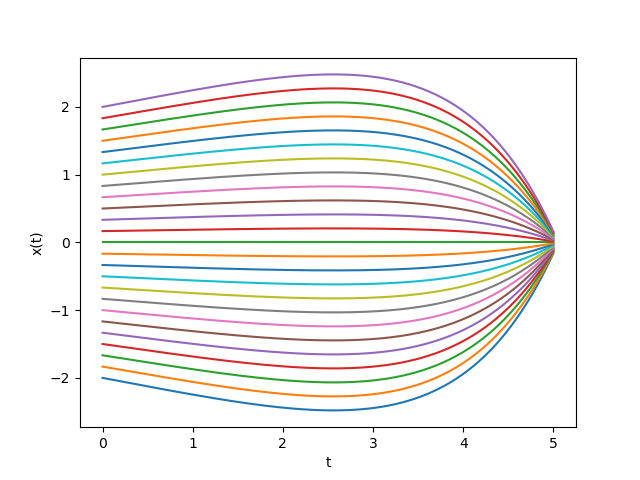}
  \centerline{\footnotesize{(d) trajectories, $T=5$}}\medskip
\end{minipage}
\caption{The open loop optimal controls and the corresponding optimal trajectories in the $16$ dimensional optimal control problem~\eqref{eqt: optctrl} with Lagrangian~\eqref{eqt:eg2_L}, source term~\eqref{eqt:eg2_f}, terminal cost~\eqref{eqt:J_min_Ji} (where $\J_i$'s are defined in~\eqref{eqt:Ji_tdep_16d}) and different terminal time $T=1$, $5$ are computed using the proposed abstract neural network architecture~\eqref{eqt: nn_u} with the implementation depicted in Fig.~\ref{fig: nn_u_RK4}. Several graphs of the first component of the optimal controls with $T=1$ and initial positions $\bx_0=(x,\bzero)\in\R^{16}$ (where $\bzero$ is the zero vector in $\R^{15}$) are shown in (a), and the first component of the corresponding optimal trajectories are shown in (b). Several graphs of the first component of the optimal controls with $T=5$ and initial positions $\bx_0=(x,\bzero)\in\R^{16}$ (where $\bzero$ is the zero vector in $\R^{15}$) are shown in (c), and the first component of the corresponding optimal trajectories are shown in (d).
\label{fig: test_tdep_16d_xu}}
\end{figure}

In what follows, we show the viscosity solution $\VNN$, the optimal controls $\optu$, and the optimal trajectories $\optx$ computed using the neural network implementations depicted in Figs.~\ref{fig: nn_V_RK4} and~\ref{fig: nn_u_RK4}. We solve a $16$ dimensional problem: $n=16$, $m=4$, and $\J$ is defined by~\eqref{eqt:J_min_Ji},
    where $\J_1,\J_2,\J_3,\J_4\colon\R^{16}\to\R$ are defined by
    \begin{equation}\label{eqt:Ji_tdep_16d}
    \begin{split}
        \J_1(\bx) &= 0.5 \|\bx\|^2 + 0.9 x_1 + 0.405,\\
        \J_2(\bx) &= 0.5 \|\bx\|^2 - 0.9 x_1 + 0.405,\\
        \J_3(\bx) &= 0.25 \|\bx\|^2 + 0.9 x_2 + 0.405,\\
        \J_4(\bx) &= 0.25 \|\bx\|^2 - 0.9 x_2 + 0.405,
    \end{split}
    \end{equation}
    for each $\bx = (x_1,\dots, x_{16})\in\R^{16}$.

The HJ PDE~\eqref{eqt: HJ_quadratic_terminal} with Hamiltonian~\eqref{eqt:eg2_H} and terminal data $\J$ defined in~\eqref{eqt:J_min_Ji} (where $\J_i$'s are defined in~\eqref{eqt:Ji_tdep_16d}) is computed using the proposed abstract neural network~\eqref{eqt: nn_S} with the implementation depicted in Fig.~\ref{fig: nn_V_RK4}. The terminal time is set to be $T=1$. The two dimensional slices for the viscosity solution are shown in Fig.~\ref{fig: test_tdep_16d}.
The subfigures (a), (b), (c), (d) show the solution at time $t=1$, $t=0.75$, $t=0.5$, $t=0.25$, respectively.

We also solve the optimal control problem with different terminal time  $T=1$ and $T=5$ using the neural network given by the abstract neural network architecture~\eqref{eqt: nn_u} with the implementation depicted in Fig.~\ref{fig: nn_u_RK4}. The open loop optimal controls $s\mapsto \optu(s|\bx_0)$ and the corresponding optimal trajectories $s\mapsto \optx(s|\bx_0)$ are then computed by solving~\eqref{eqt: ODE_optx} with the fourth order Runge-Kutta method whose one-step neural network representation is shown in Fig.~\ref{fig: nn_RK4}.
The graphs of the first component of the optimal controls $\optu$ and the first component of the optimal trajectories $\optx$ to the $16$ dimensional problem with the initial position $\bx_0=(x,\bzero)\in\R^{16}$ (where $\bzero$ is the zero vector in $\R^{15}$) and the terminal cost~\eqref{eqt:J_min_Ji} (where $\J_i$'s are defined in~\eqref{eqt:Ji_tdep_16d}) are shown in Fig.~\ref{fig: test_tdep_16d_xu}. In the figure, (a) and (c) show the optimal controls $\optu$ with $T=1$ and $T=5$, respectively, while (b) and (d) show the corresponding optimal trajectories $\optx$ with $T=1$ and $T=5$, respectively.

\subsection{An optimal control problem in Newton mechanics} \label{subsec: test_newton}
We consider the optimal control problem~\eqref{eqt: optctrl} whose Lagrangian $L$ reads
\begin{equation}\label{eqt:eg3_L}
    L(t,\bx,\bu) = \frac{1}{2} \|\bx - \bx_r(t)\|^2 + \frac{1}{2000} \|\bu\|^2 \quad \forall t\in[0,T], \bx\in\Rn, \bu\in\R^l,
\end{equation}
where we set $n=2l$ and define the function $\bx_r\colon [0,T]\to \Rn$ by
\begin{equation*}
    \bx_r(t) = 5\sin t \begin{pmatrix}
        I_l \\
        \Ol
    \end{pmatrix} + 5\cos t \begin{pmatrix}
        \Ol \\
        I_l
    \end{pmatrix}.
\end{equation*}
The source term $f$ in~\eqref{eqt: ode} is defined by
\begin{equation}\label{eqt:eg3_f}
    f(t,\bx,\bu) = \begin{pmatrix}
        \Ol & I_l \\
        \Ol & \Ol
    \end{pmatrix}\bx+ \begin{pmatrix}
        \Ol \\
        I_l
    \end{pmatrix}\bu \quad \forall t\in[0,T], \bx\in\Rn, \bu\in\R^l.
\end{equation}
If we denote $\bx(s)=(\bx_1(s),\bx_2(s))$, where $\bx_1(s),\bx_2(s)\in\R^l$ for each $s\in[t_0,T]$, then the Cauchy problem~\eqref{eqt: ode} becomes
\begin{equation*}
\begin{dcases}
\dot{\bx}_1(s) = \bx_2(s) & s\in (\tzero,T),\\
\dot{\bx}_2(s) = \bu(s) &s\in (\tzero,T),\\
(\bx_1(t_0), \bx_2(t_0)) = \bx_0.
\end{dcases}
\end{equation*}
This is the ODE in Newton mechanics, where $\bx_1$ denotes the position of a particle, $\bx_2$ denotes its velocity, and $\bu$ denotes its acceleration.
The corresponding HJ PDE is in the form of~\eqref{eqt: HJ_quadratic_terminal} where the Hamiltonian $H$ is defined by
\begin{equation}\label{eqt:eg3_H}
    H(t,\bx,\bp) = \frac{1}{2}\langle \bp, \Cpp \bp\rangle  - \langle \bp, \Cxp \bx\rangle - \frac{1}{2}\|\bx - \bx_r(t)\|^2
    \quad \forall t\in[0,T], \bx,\bp\in\Rn,
\end{equation}
where the coefficients $\Cpp,\Cxp\in\matnn$ are constant matrices given by
\begin{equation*}
    \Cpp = 1000\begin{pmatrix}
        \Ol & \Ol \\
        \Ol & I_l
    \end{pmatrix}
    \quad\text{ and }\quad
    \Cxp = \begin{pmatrix}
        \Ol & I_l \\
        \Ol & \Ol
    \end{pmatrix}.
\end{equation*}
We consider the terminal data $\J\colon \Rn\to\R$ defined by
\begin{equation} \label{eqt: newton_J}
    \J(\bx) = \min\left\{\frac{1}{320}\left((x_1 + 2)^2 + \sum_{i=2}^nx_i^2\right), \frac{1}{320}\left((x_1 - 2)^2 + \sum_{i=2}^nx_i^2\right)\right\},
\end{equation}
for each $\bx=(x_1,x_2,\dots, x_n)\in\Rn$.

Note that this problem requires a slight modification of the HJ PDE~\eqref{eqt: HJ_quadratic_terminal} and the optimal control problem~\eqref{eqt: optctrl_quadratic} considered in this paper, because of the term $\frac{1}{2}\|\bx-\bx_r(t)\|^2$.
However, the two abstract neural network architectures~\eqref{eqt: nn_S} and~\eqref{eqt: nn_u} can still be used to compute the viscosity solution and the optimal control, where in the $i$-th neuron, the function $P_i\in C([0,T];\symn)$ solves the Riccati FVP~\eqref{eqt: odeP} which reads
\begin{equation*}
    \left\{
    \begin{aligned}
    &\dot{\Sxx}_i(t) =  1000\Sxx_i(t)^T\begin{pmatrix}
        \Ol & \Ol \\
        \Ol & I_l
    \end{pmatrix}\Sxx_i(t) - \Sxx_i(t)^T\begin{pmatrix}
        \Ol & I_l \\
        \Ol & \Ol
    \end{pmatrix} \\
    &\quad\quad\quad\quad\quad\quad - \begin{pmatrix}
        \Ol & \Ol \\
        I_l & \Ol
    \end{pmatrix}\Sxx_i(t) - I_n &t\in(0,T),\\
    &\Sxx_i(T) = \frac{1}{160}I_n,
    \end{aligned}
    \right.
\end{equation*}
the functions $\bq_i\in C(0,T;\Rn)$ solves the modified FVP which reads
\begin{equation*}
    \begin{dcases}
    \dot{\Sx}_i(t) = 1000\Sxx_i(t)^T\begin{pmatrix}
        \Ol & \Ol \\
        \Ol & I_l
    \end{pmatrix}\Sx_i(t) - \begin{pmatrix}
        \Ol & \Ol \\
        I_l & \Ol
    \end{pmatrix}\Sx_i(t) + \bx_r(t) &t\in(0,T),\\
    \Sx_i(T) = \Gx_i,
    \end{dcases}
\end{equation*}
and $r_i\in C(0,T;\R)$ solves the modified FVP which reads
\begin{equation*}
    \begin{dcases}
    \dot{\Sc}_i(t) = 500\Sx_i(t)^T\begin{pmatrix}
        \Ol & \Ol \\
        \Ol & I_l
    \end{pmatrix}\Sx_i(t) - \frac{25l}{2} &t\in(0,T),\\
    \Sc_i(T) = \Gc_i.
    \end{dcases}
\end{equation*}
For the specific terminal data $\J$ in~\eqref{eqt: newton_J}, $\Gx_i\in\Rn$ and $\Gc_i\in\R$ are given by
\begin{equation*}
    \Gx_1 = \frac{1}{80}(1,\bzero)^T, \quad \Gc_1 = \frac{1}{80},
    \quad \Gx_2 = -\frac{1}{80}(1,\bzero)^T, \quad \Gc_2 = \frac{1}{80},
\end{equation*}
where $\bzero$ denotes the zero vector in $\R^{n-1}$.

\begin{figure}[htbp]
\begin{minipage}[b]{.49\linewidth}
  \centering
  \includegraphics[width=.99\textwidth]{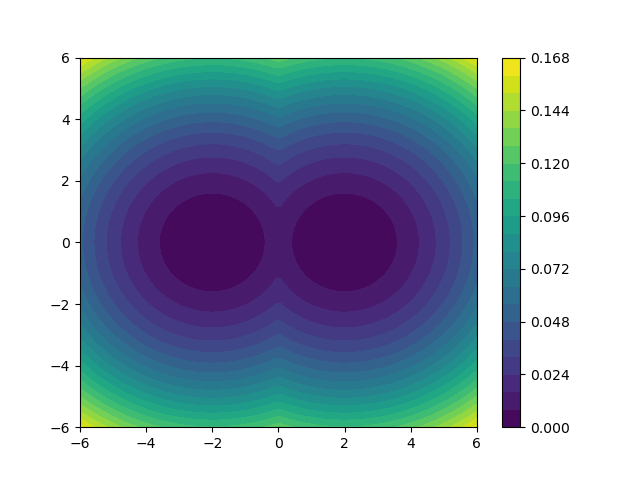}
  \centerline{\footnotesize{(a) $t=1$}}\medskip
\end{minipage}
\hfill  
\begin{minipage}[b]{.49\linewidth}
  \centering
  \includegraphics[width=.99\textwidth]{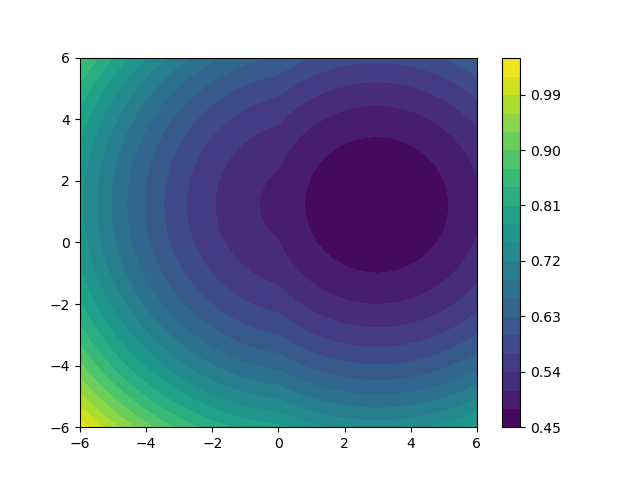}
  \centerline{\footnotesize{(b) $t=0.995$}}\medskip
\end{minipage}
\begin{minipage}[b]{.49\linewidth}
  \centering
    \includegraphics[width=.99\textwidth]{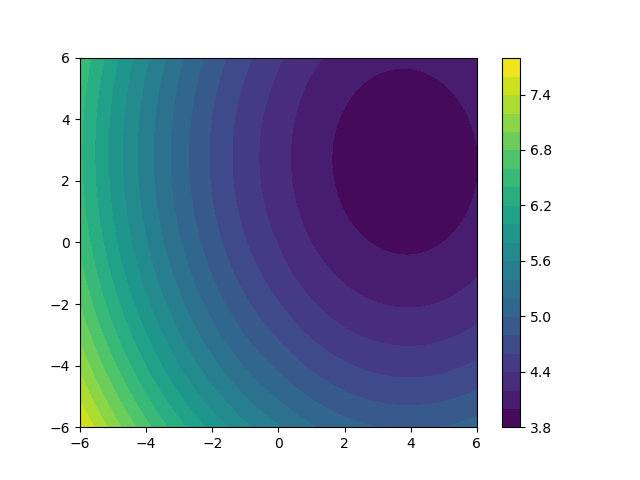}
  \centerline{\footnotesize{(c) $t=0.95$}}\medskip
\end{minipage}
\hfill
\begin{minipage}[b]{.49\linewidth}
  \centering
    \includegraphics[width=.99\textwidth]{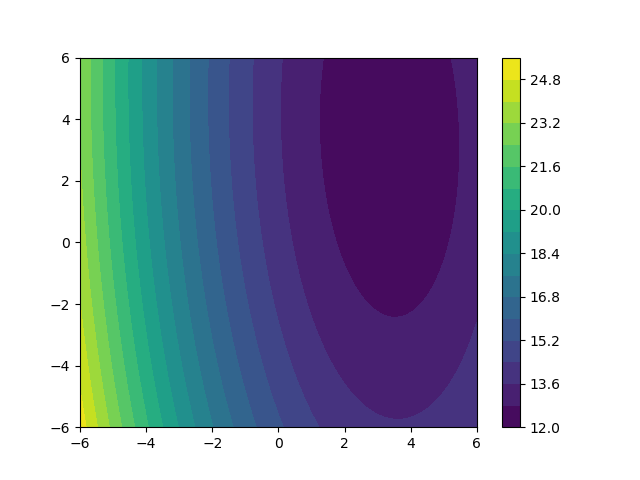}
  \centerline{\footnotesize{(d) $t=0.75$}}\medskip
\end{minipage}
\caption{The viscosity solution $\VNN$ to the $16$ dimensional HJ PDE~\eqref{eqt: HJ_quadratic_terminal} with Hamiltonian~\eqref{eqt:eg3_H}, terminal data~\eqref{eqt: newton_J} and terminal time $T=1$ is computed using the proposed abstract neural network architecture~\eqref{eqt: nn_S} (depicted in Fig.~\ref{fig: nn_Riccati}) with the implementation depicted in Fig.~\ref{fig: nn_V_RK4}. 
The two dimensional slices of $\VNN$ at time $t=1$ (i.e., the terminal cost), $t=0.995$, $t=0.95$ and $t=0.75$ 
are shown in the subfigures (a), (b), (c) and (d), respectively.
The color in each subfigure shows the solution value $\VNN(t, \bx)$, where the spatial variable $\bx$ is in the form of $(x_1,\mathbf{0},x_2,\mathbf{0})\in\R^{16}$ (with $\mathbf{0}$ denoting the zero vector in $\R^7$) for some points $x_1\in\R$ and $x_2\in\R$ which are represented by the $x$ and $y$ axes.
\label{fig: test_newton}}
\end{figure}

Here, we show the numerical results for $l=8$ and $n=2l=16$. The viscosity solution $\VNN$ with terminal time  $T=1$ is computed using the abstract neural network architecture~\eqref{eqt: nn_S} (depicted in Fig.~\ref{fig: nn_Riccati}) with the implementation shown in Fig.~\ref{fig: nn_V_RK4}. The two dimensional slices of $\VNN$ at $t = 1$, $0.995$, $0.95$, $0.75$
are plotted in Fig.~\ref{fig: test_newton} (a), (b), (c), (d), respectively. 

\begin{figure}[htbp]
\begin{minipage}[b]{.49\linewidth}
  \centering
  \includegraphics[width=.99\textwidth]{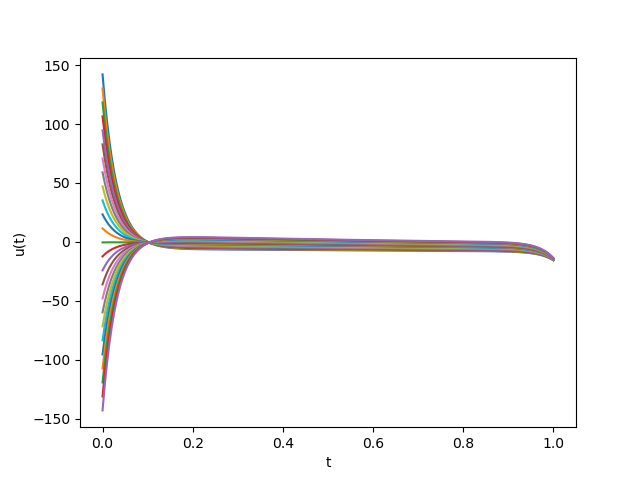}
  \centerline{\footnotesize{(a) $T=1$}}\medskip
\end{minipage}
\hfill  
\begin{minipage}[b]{.49\linewidth}
  \centering
  \includegraphics[width=.99\textwidth]{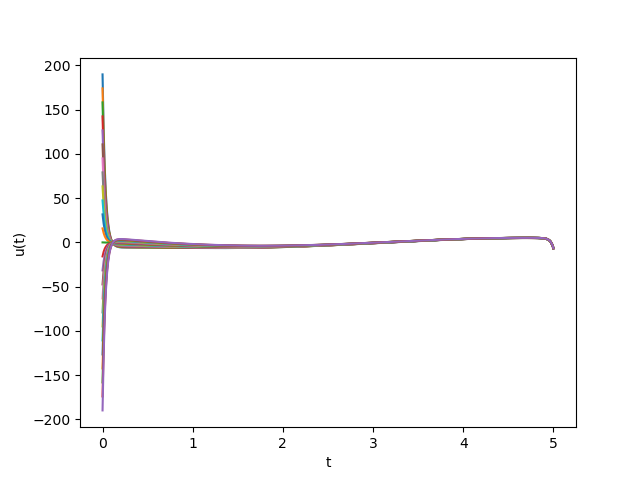}
  \centerline{\footnotesize{(b) $T=5$}}\medskip
\end{minipage}
\begin{minipage}[b]{.49\linewidth}
  \centering
    \includegraphics[width=.99\textwidth]{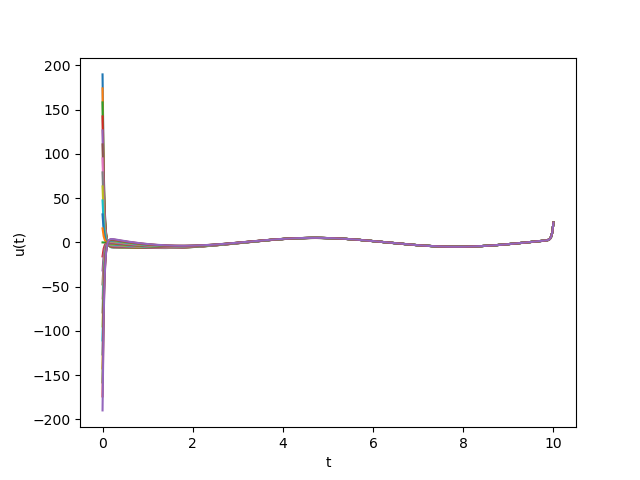}
  \centerline{\footnotesize{(c) $T=10$}}\medskip
\end{minipage}
\caption{The open loop optimal controls in the $16$ dimensional optimal control problem~\eqref{eqt: optctrl} with Lagrangian~\eqref{eqt:eg3_L}, source term~\eqref{eqt:eg3_f}, terminal cost~\eqref{eqt: newton_J} and different terminal time $T=1,5,10$ are computed using the proposed abstract neural network architecture~\eqref{eqt: nn_u} depicted in Fig.~\ref{fig: nn_u} with the implementation shown in Fig.~\ref{fig: nn_u_RK4}. Several graphs of the first component of the optimal controls with $T=1$, $T=5$ and $T=10$ are shown in (a), (b) and (c), respectively. In each figure, different trajectories correspond to different initial positions $\bx_0=(x,\bzero)\in\R^{16}$ where $\bzero$ is the zero vector in $\R^{15}$.
\label{fig: test_newton_u}}
\end{figure}

\begin{figure}
    \begin{minipage}[b]{.49\linewidth}
  \centering
    \includegraphics[width=.99\textwidth]{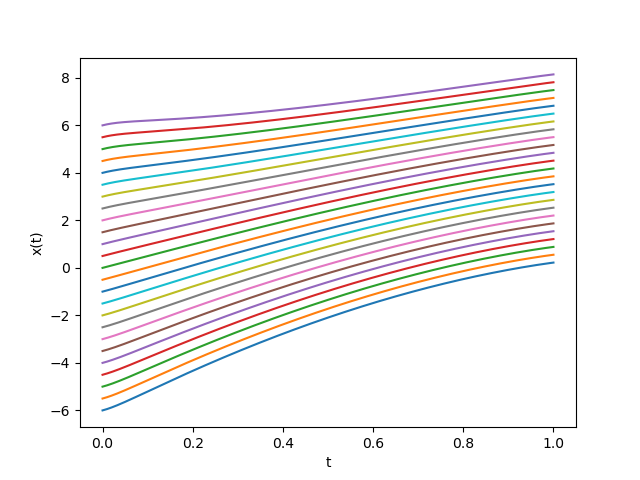}
  \centerline{\footnotesize{(a) 1st component $x_1^*$, $T=1$}}\medskip
\end{minipage}
\hfill
\begin{minipage}[b]{.49\linewidth}
  \centering
    \includegraphics[width=.99\textwidth]{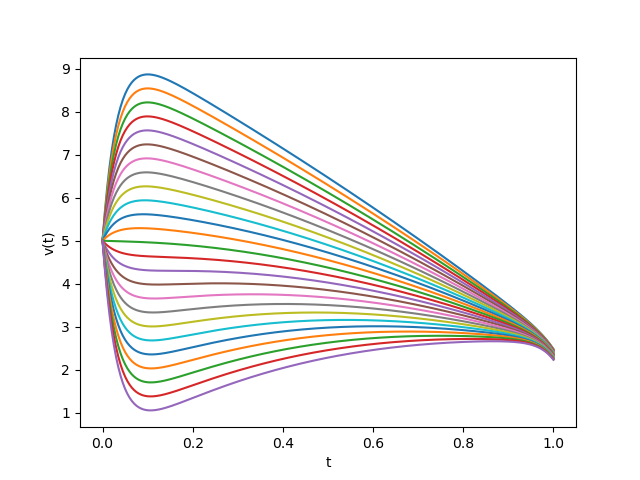}
  \centerline{\footnotesize{(b) 9th component $x_9^*$, $T=1$}}\medskip
\end{minipage}
\begin{minipage}[b]{.49\linewidth}
  \centering
    \includegraphics[width=.99\textwidth]{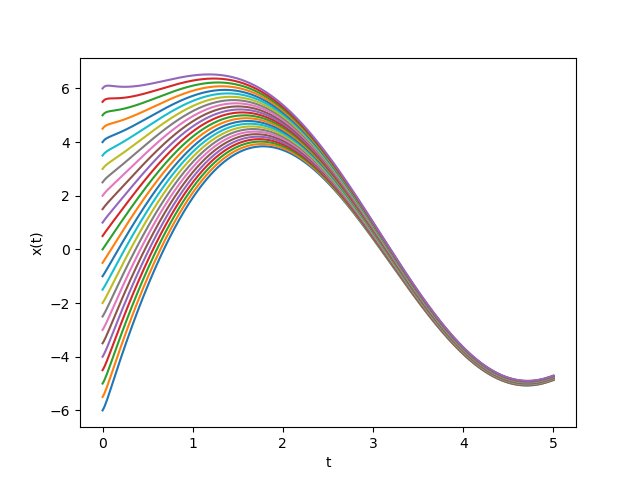}
  \centerline{\footnotesize{(c) 1st component $x_1^*$, $T=5$}}\medskip
\end{minipage}
\hfill
\begin{minipage}[b]{.49\linewidth}
  \centering
    \includegraphics[width=.99\textwidth]{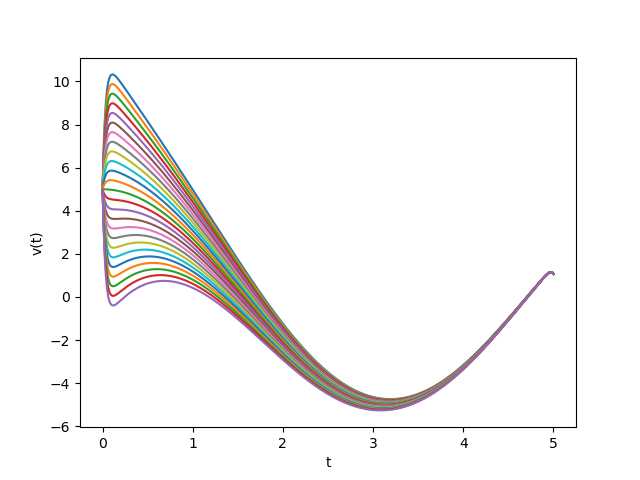}
  \centerline{\footnotesize{(d) 9th component $x_9^*$, $T=5$}}\medskip
\end{minipage}
\begin{minipage}[b]{.49\linewidth}
  \centering
    \includegraphics[width=.99\textwidth]{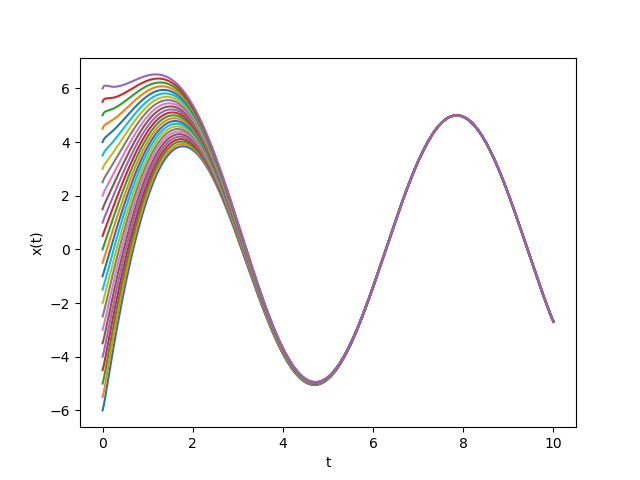}
  \centerline{\footnotesize{(e) 1st component $x_1^*$, $T=10$}}\medskip
\end{minipage}
\hfill
\begin{minipage}[b]{.49\linewidth}
  \centering
    \includegraphics[width=.99\textwidth]{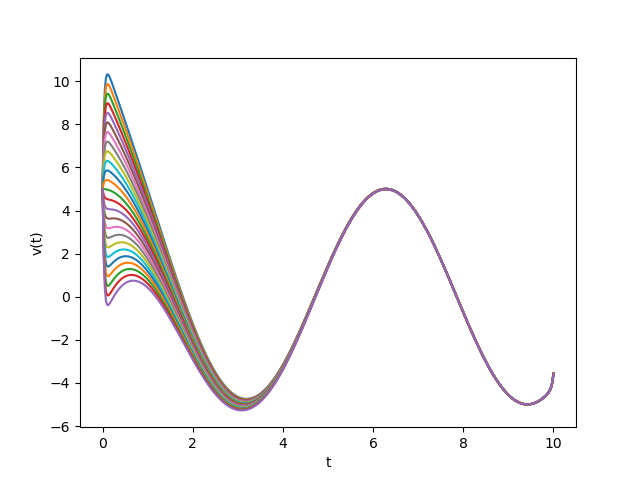}
  \centerline{\footnotesize{(f) 9th component $x_9^*$, $T=10$}}\medskip
\end{minipage}
    \caption{The corresponding optimal trajectories $\optx$ in the $16$ dimensional optimal control problem~\eqref{eqt: optctrl} with Lagrangian~\eqref{eqt:eg3_L}, source term~\eqref{eqt:eg3_f}, terminal cost~\eqref{eqt: newton_J} and different terminal time $T=1,5,10$ are computed using the proposed abstract neural network architecture~\eqref{eqt: nn_u} depicted in Fig.~\ref{fig: nn_u} with the implementation shown in Fig.~\ref{fig: nn_u_RK4}. Several graphs of the first component $x_1^*$ of the optimal trajectories $\optx$ with $T=1$, $T=5$ and $T=10$ are shown in (a), (c) and (e), respectively. The corresponding graphs of the ninth component $x_9^*$ of the optimal trajectories $\optx$ with $T=1$, $T=5$ and $T=10$ are shown in (b), (d) and (f), respectively. In each figure, different trajectories correspond to different initial positions $\bx_0=(x,\bzero)\in\R^{16}$ where $\bzero$ is the zero vector in $\R^{15}$.
    \label{fig: test_newton_xv}}
\end{figure}

The optimal controls with different terminal time  $T=1,5,10$ are computed using the abstract neural network architecture~\eqref{eqt: nn_u} depicted in Fig.~\ref{fig: nn_u} with the implementation depicted in Fig.~\ref{fig: nn_u_RK4}. 
The open loop optimal controls and the corresponding optimal trajectories are computed by solving~\eqref{eqt: ODE_optx} with the fourth order Runge-Kutta method whose one-step neural network representation is shown in Fig.~\ref{fig: nn_RK4}.
The graphs of the first components of the optimal trajectories with terminal time $T=1$, $T=5$, $T=10$ and different initial positions $\bx_0=(x,\bzero)\in\R^{16}$ (where $\bzero$ is the zero vector in $\R^{15}$) are shown in Fig.~\ref{fig: test_newton_u} (a), (b), (c), respectively. The graphs of the first components of the corresponding optimal trajectories with $T=1, 5, 10$ are shown in Fig.~\ref{fig: test_newton_xv} (a), (c), (e), while the graphs of the ninth components of the optimal trajectories with $T=1, 5, 10$ are shown in Fig.~\ref{fig: test_newton_xv} (b), (d), (f).
From the optimal controls and trajectories for longer time horizons, it appears there is a turnpike phenomenon.

\subsection{An optimal control problem with general terminal cost}
\label{subsec:example_admm}

In this section, we consider a more general terminal cost $\J$ in the form of~\eqref{eqt:J_min_Ji} with $m=2$,
\begin{equation}\label{eqt:eg4_defJ1_J2}
    \J_1(\bx) = \|\bx - \bx_1\|_1,\quad 
    \J_2(\bx) = \|\bx - \bx_2\|,
\end{equation}
for all $\bx\in\Rn$. Recall that $\|\cdot\|_1$ and $\|\cdot\|$ denote the $\ell^1$-norm and $\ell^2$-norm in $\Rn$, respectively.
For illustration purposes, we set $\bx_1 = (1,1,\mathbf{0})$ and $\bx_2 = (-1,-1,\mathbf{0})$. 
We solve the same optimal control problem as in Section~\ref{subsec: test_tdep}, with the terminal cost~\eqref{eqt:J_min_Ji} (where $\J_i$'s are defined in~\eqref{eqt:eg4_defJ1_J2}).
The value function $\VNN$ and the optimal control $\bu$ are computed using the neural network architecture in Figs.~\ref{fig: nn_V_RK4} and~\ref{fig: nn_u_RK4}, respectively, where the parameters $\Gxx_i, \Gx_i$ and $\Gc_i$ are trained using the ADMM algorithm described in Section~\ref{sec:admm_method}.
Note that in each iteration in the ADMM algorithm, we need to solve the optimal control problem~\eqref{eqt:admm_step1} once using our proposed neural network architecture in Fig.~\ref{fig: nn_V_RK4}.

\begin{figure}[htbp]
\begin{minipage}[b]{.49\linewidth}
  \centering
  \includegraphics[width=.99\textwidth]{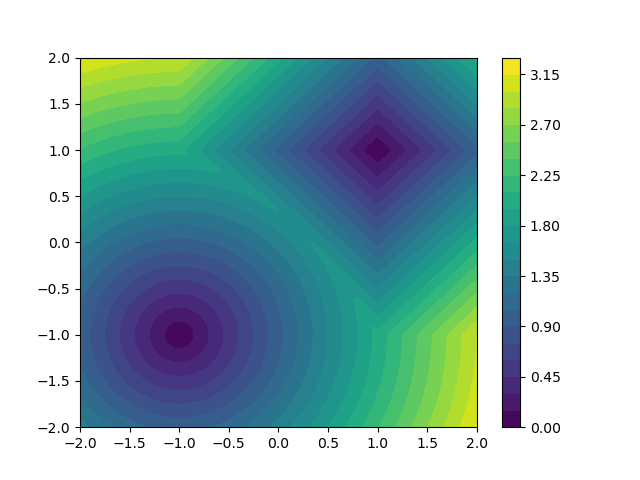}
  \centerline{\footnotesize{(a) $t=1$}}\medskip
\end{minipage}
\hfill  
\begin{minipage}[b]{.49\linewidth}
  \centering
  \includegraphics[width=.99\textwidth]{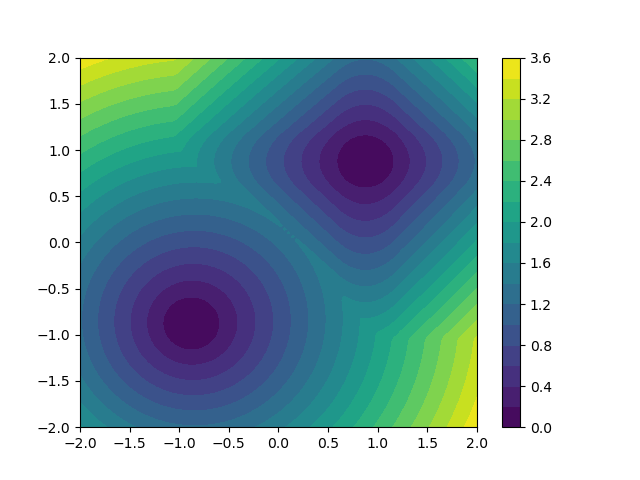}
  \centerline{\footnotesize{(b) $t=0.75$}}\medskip
\end{minipage}
\begin{minipage}[b]{.49\linewidth}
  \centering
    \includegraphics[width=.99\textwidth]{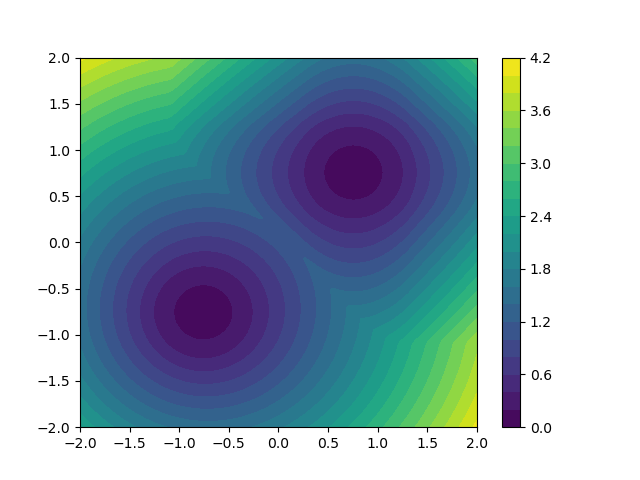}
  \centerline{\footnotesize{(c) $t=0.5$}}\medskip
\end{minipage}
\hfill
\begin{minipage}[b]{.49\linewidth}
  \centering
    \includegraphics[width=.99\textwidth]{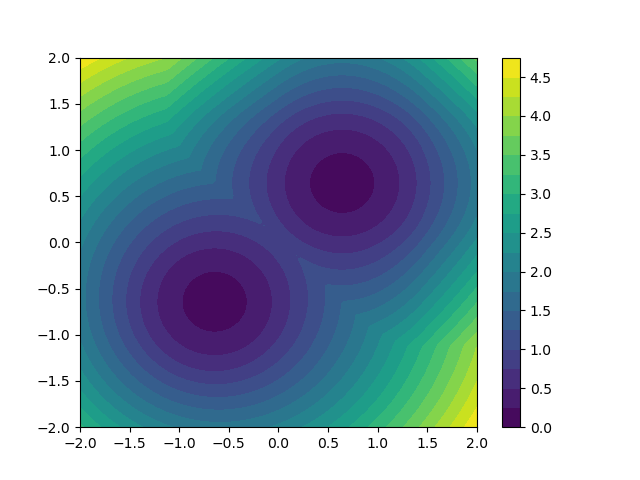}
  \centerline{\footnotesize{(d) $t=0.25$}}\medskip
\end{minipage}
\caption{The viscosity solution $\VNN$ to the $16$ dimensional HJ PDE with Hamiltonian~\eqref{eqt:eg2_H}, terminal data~\eqref{eqt:J_min_Ji} (where $\J_i$'s are defined in~\eqref{eqt:eg4_defJ1_J2}) and terminal time $T=1$ is computed using the proposed abstract neural network architecture depicted in Fig.~\ref{fig: nn_V_RK4}, whose parameters are trained using the ADMM method. 
The two dimensional slices of $\VNN$ at time $t=1$ (i.e., the terminal cost), $t=0.75$, $t=0.5$ and $t=0.25$ are shown in the subfigures (a), (b), (c) and (d), respectively.
The color in each subfigure shows the solution value $\VNN(t, \bx)$, where the spatial variable $\bx$ is in the form of $(x_1,x_2,\bzero)\in\R^{16}$ (with $\bzero$ denoting the zero vector in $\R^{14}$) for some points $x_1\in\R$ and $x_2\in\R$ which are represented by the $x$ and $y$ axes.
\label{fig:eg4_V}}
\end{figure}

\begin{figure}[htbp]
\begin{minipage}[b]{.49\linewidth}
  \centering
  \includegraphics[width=.99\textwidth]{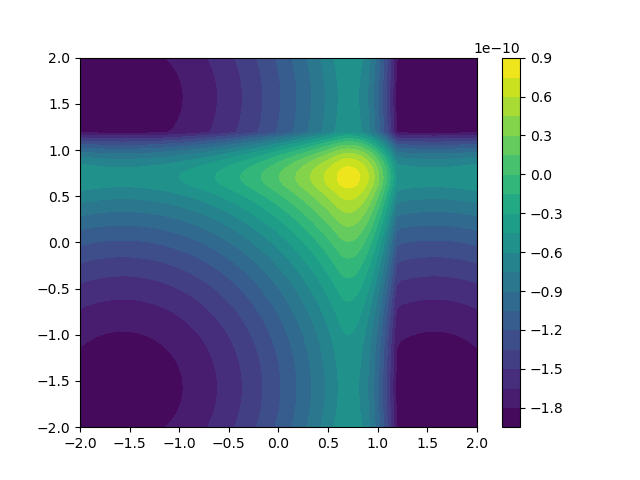}
  \centerline{\footnotesize{(a)}}\medskip
\end{minipage}
\hfill  
\begin{minipage}[b]{.49\linewidth}
  \centering
  \includegraphics[width=.99\textwidth]{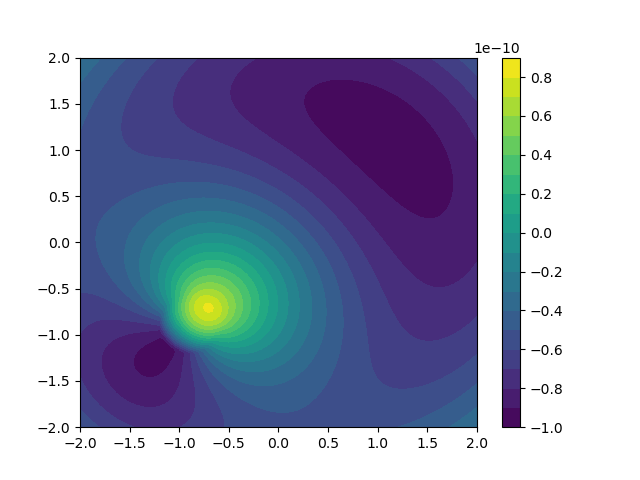}
  \centerline{\footnotesize{(b)}}\medskip
\end{minipage}
\begin{minipage}[b]{.49\linewidth}
  \centering
    \includegraphics[width=.99\textwidth]{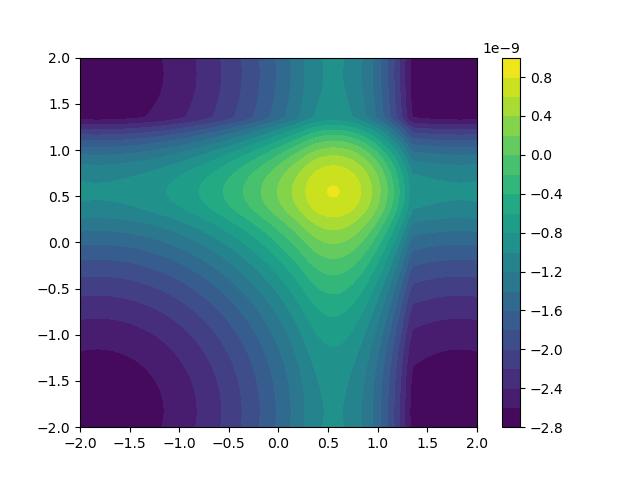}
  \centerline{\footnotesize{(c)}}\medskip
\end{minipage}
\hfill
\begin{minipage}[b]{.49\linewidth}
  \centering
    \includegraphics[width=.99\textwidth]{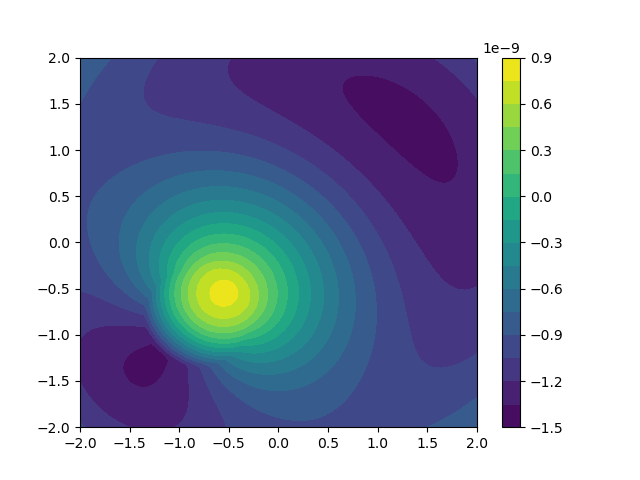}
  \centerline{\footnotesize{(d)}}\medskip
\end{minipage}
\begin{minipage}[b]{.49\linewidth}
  \centering
    \includegraphics[width=.99\textwidth]{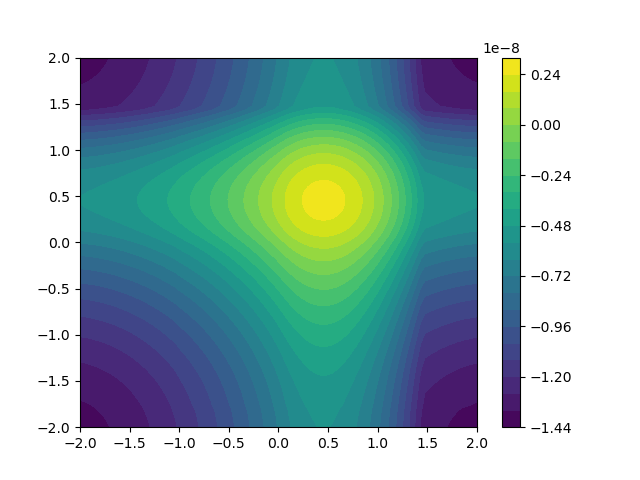}
  \centerline{\footnotesize{(e)}}\medskip
\end{minipage}
\hfill
\begin{minipage}[b]{.49\linewidth}
  \centering
    \includegraphics[width=.99\textwidth]{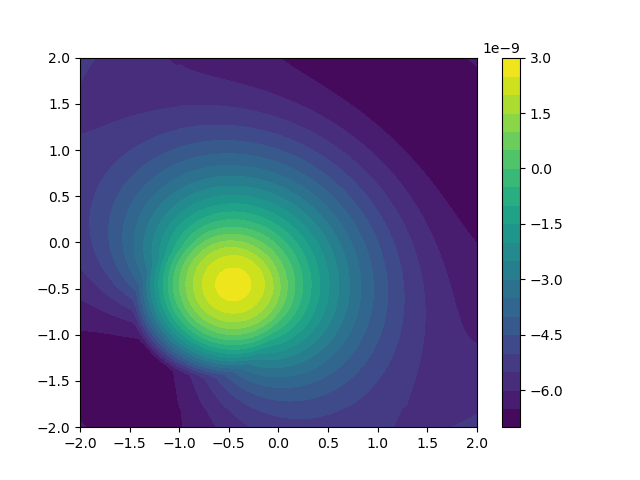}
  \centerline{\footnotesize{(f)}}\medskip
\end{minipage}
\caption{The residual $-\frac{\partial \V_i}{\partial t} + H(t,\bx,\nabla_{\bx}\V_i)$ in the HJ PDE with terminal time $T=1$, Hamiltonian~\eqref{eqt:eg2_H} and terminal condition~\eqref{eqt:J_min_Ji} (where $\J_i$'s are defined in~\eqref{eqt:eg4_defJ1_J2}) is shown for each $i\in\{1,2\}$ at different time $t$, where $\V_i$ is the solution to the $i$-th subproblem. Figures (a), (c), (e) show the residuals for $\V_1$ at time $t=0.75T$, $t=0.5T$ and $t=0.25T$, while figures (b), (d), (f) show the residuals for $\V_2$ at time $t=0.75T$, $t=0.5T$ and $t=0.25T$, respectively. 
In each subfigure, we show the two dimensional slices of the residual function. The color shows the residual value at $(t, \bx)$, where the spatial variable $\bx$ is in the form of $(x_1,x_2,\bzero)\in\R^{16}$ (with $\bzero$ denoting the zero vector in $\R^{14}$) for some points $x_1\in\R$ and $x_2\in\R$ which are represented by the $x$ and $y$ axes.
\label{fig:eg4_pde_err}}
\end{figure}

We show the numerical results for the $16$-dimensional problem. The two dimensional slices of the viscosity solution $\VNN$ with terminal time $T=1$ at time $t=1$, $0.75$, $0.5$, $0.25$ are plotted in Fig.~\ref{fig:eg4_V} (a), (b), (c), (d), respectively. 
The residual $-\frac{\partial \V_i}{\partial t} + H(t,\bx,\nabla_{\bx}\V_i)$ of the HJ PDE with terminal condition $\J_1$ and $\J_2$ at different time $t=0.75$, $0.5$, $0.25$ are shown in Fig.~\ref{fig:eg4_pde_err}. We can observe a small error from these error plots, which numerically validate that our proposed neural network architecture indeed solves the viscosity solution to the corresponding HJ PDE.

\begin{figure}[htbp]
\begin{minipage}[b]{.49\linewidth}
  \centering
  \includegraphics[width=.99\textwidth]{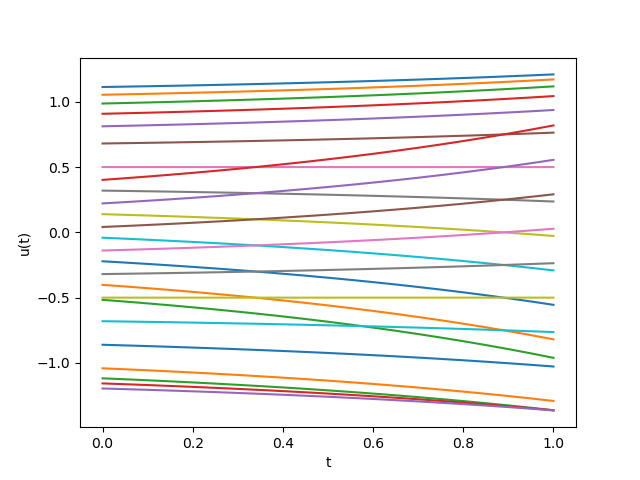}
  \centerline{\footnotesize{(a)}}\medskip
\end{minipage}
\hfill  
\begin{minipage}[b]{.49\linewidth}
  \centering
  \includegraphics[width=.99\textwidth]{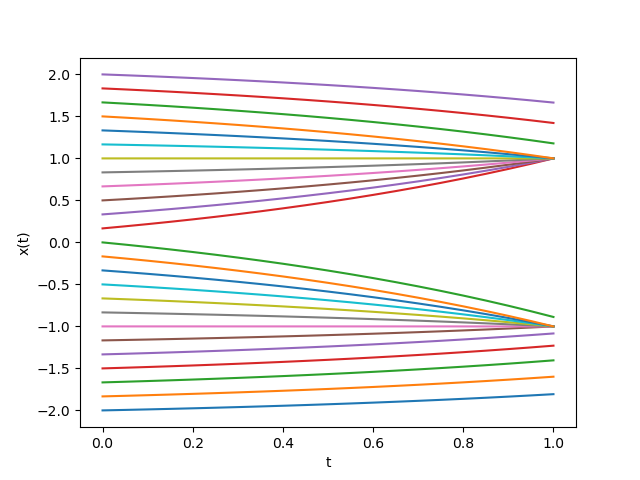}
  \centerline{\footnotesize{(b)}}\medskip
\end{minipage}
\begin{minipage}[b]{.49\linewidth}
  \centering
    \includegraphics[width=.99\textwidth]{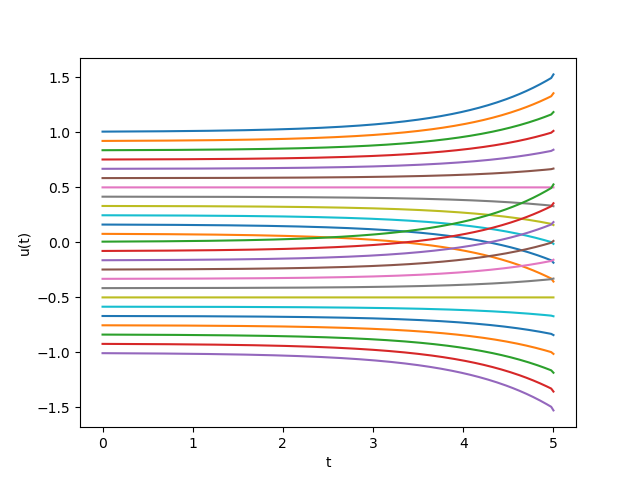}
  \centerline{\footnotesize{(c)}}\medskip
\end{minipage}
\hfill
\begin{minipage}[b]{.49\linewidth}
  \centering
    \includegraphics[width=.99\textwidth]{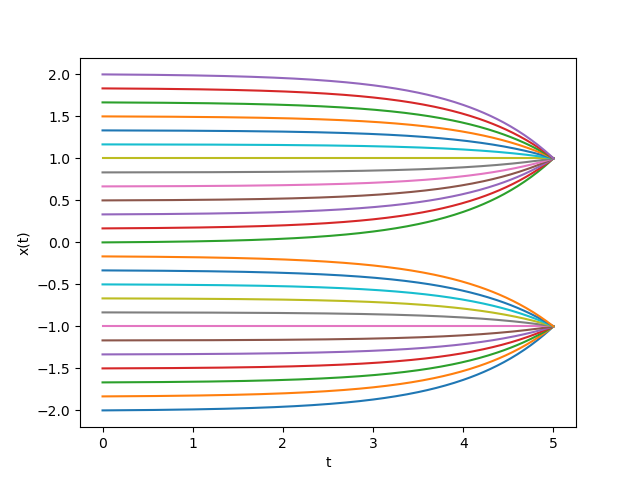}
  \centerline{\footnotesize{(d)}}\medskip
\end{minipage}
\caption{The open loop optimal controls and the optimal trajectories in the $16$ dimensional optimal control problem with Lagrangian~\eqref{eqt:eg2_L}, source term~\eqref{eqt:eg2_f}, terminal cost~\eqref{eqt:J_min_Ji} (where $\J_i$'s are defined in~\eqref{eqt:eg4_defJ1_J2}) and different terminal time $T=1,5$ are computed using the proposed abstract neural network architecture~\eqref{eqt: nn_u} with the implementation depicted in Fig.~\ref{fig: nn_u_RK4}. Several graphs of the optimal controls with $T=1$ are shown in (a), whose corresponding optimal trajectories are shown in (b). Several graphs of the optimal controls with $T=5$ are shown in (c), whose corresponding optimal trajectories are shown in (d).
\label{fig:eg4_xu}}
\end{figure}

We also compute several optimal controls and trajectories with different initial position $\bx=(x, \bzero)\in\R^n$ and fixed initial time $t_0 = 0$, and the graphs of their first components are shown in Fig.~\ref{fig:eg4_xu}. The optimal controls with terminal time $T=1$, $5$ are shown in Fig.~\ref{fig:eg4_xu} (a), (c), while the optimal trajectories with terminal time $T=1$, $5$ are shown in Fig.~\ref{fig:eg4_xu} (b), (d), respectively.

 \subsection{An FPGA implementation and numerical results}
 \label{subsec:FPGA}

We now briefly describe am implementation of our proposed neural network on FPGA to illustrate the performance that can be achieved using simple precision floating points. Specifically, we only present an FPGA implementation with low latency, where the latency corresponds to the amount of time the neural network takes to produce one result.

FPGAs are an array of programmable logic blocks and memory elements that are connected together using a programmable interconnect. FPGAs contain different types of logic resources. These resources include general purposes logics such as lookup tables (LUTs) and Flip-Flops (FFs), more specialized arithmetic units, such as digital signal processing units (DSPs), and memory such as Block Random Access Memory (BRAMs). We refer the reader to~\cite{2018arXiv180503648K} for a concise description of FPGAs. We use the Xilinx Alveo U280 board with a target design running at 300 MHz.
The main computational burden of our proposed neural network consists of matrix-matrix multiplications used in the fourth order Runge-Kutta method for solving an FVP as described in Fig.~\ref{fig: nn_RK4}. Traditional non-parallel algorithms for performing matrix-matrix multiplications have an $O(n^3)$ time complexity. Using the parallel programming feature of FPGAs we can obtain a complexity of $O(n^2)$ for computing matrix-matrix multiplications (see~\cite{2018arXiv180503648K} for instance). Therefore, we spend most of FPGA resources on performing these matrix-matrix multiplications in order to reduce the latency of the design. Note that the Alveo U280 board is composed of three ``chiplets" and crossing chiplets consumes scarce routing resources that severely degrades performance and prevents scaling. Therefore, we only consider FPGA designs that use less that 30\% of the FPGA resources so that no chiplet is crossed. Table~\ref{tab:FPGAfloat} presents the FPGA resources and latencies to implement our proposed neural network depicted in Fig.~\ref{fig: nn_V_RK4} for various dimensions $n$ and numbers of layers $L$. We observe from the table that the latencies for $(n,L)= (16,8)$, $(32,4)$, and $(64,2)$ are  2.1110e-05s, 7.5150e-05s, and 2.8600e-04s, respectively. We also implemented our proposed neural network architecture on CPUs using C++ to highlight the boost of performance we can obtain using FPGAs. We perform 1,000,000 runs on a single Intel core I7-1165G7 and report the average time to produce a result for $(n,L)= (16,8)$, $(32,4)$, $(64,2)$ in Table~\ref{tab:FPGAtime} as well as the speed-up compared to our FPGA implementation. We observe a speed-up from 12 to about 20 depending on the dimension $n$
and the number of layers $L$.
Our FPGA design also allows for larger number of layers than those reported here. We simply iterate the FPGA kernel that we designed here for the neural network with fewer layers. 
In these cases, the amount of FPGA resources remain the same but the latency is multiplied by the number of iterations of the FPGA kernel. 
\begin{table}[]
    \centering
 \begin{tabular}{c|c|c|c|c|c}
 \hline
$n/L$ &  Latency (ns)& BRAMs      & DSPs        &   FFs           & LUTs\\
\hline
\hline
16/8& 6,345 (2.111E4)  & 601(14\%)   & 2402(26\%) & 354,013(13\%)&258,710(19\%) \\  
32/4& 22,547(7.515E4)  & 602(14\%) & 2482(27\%) & 353,225(13\%)&248,369(19\%)  \\ 
64/2& 85,547 (2.860E5) & 608(15\%) & 2522(27\%) & 352,715(13\%)&242,886(19\%)  \\ 
\hline
 \end{tabular}
    \caption{FPGA resources and latencies in cycles and nanoseconds (ns) to implement $L$ layers of the neural network for various dimensions $n$  using simple precision floating point on a Xilinx Alveo U280 board with a frequency of 300 MHz.}
    \label{tab:FPGAfloat}
\end{table}

\begin{table}[]
    \centering
 \begin{tabular}{c|c|c|c}
 \hline
n/L       &          CPU time  &        FPGA  time  &  speed up\\
\hline
\hline
16/8    & 2.6310e-04s&  2.1110e-05s&    12.463\\ 
32/4   & 1.2021e-03s&  7.5150e-05s &   15.996\\
64/2    & 5.9730e-03s&  2.8600e-04s &   20.885 \\
\hline
 \end{tabular}
    \caption{Comparison of the average time for 1,000,000 runs for various dimensions and number of layers on a single Intel Core I7-1165G7 and our FPGA implementation on a Xilinx Alveo U280 board running at 300 MHz. The speed-up using FPGA compared to the CPU is presented in the last column.}
    \label{tab:FPGAtime}
\end{table}

\section{Conclusion}\label{sec:conclusion}
We propose two abstract neural network architectures depicted in Figs.~\ref{fig: nn_Riccati} and~\ref{fig: nn_u}, which respectively solve certain high dimensional HJ PDEs and are used to compute the optimal controls in the corresponding optimal control problems. To implement these abstract architectures, we present two Resnet-type deep neural network implementations and show several numerical results in Section~\ref{sec:implementation}. These architectures pave the way to leverage dedicated hardware designed for neural networks to obtain efficient implementations of the numerical algorithms for certain optimal control problems and HJ PDEs. It has potential in real-time computations for these high dimensional problems. 
Moreover, these architectures are designed based on the theories of linear-quadratic controls and min-plus algebra, and hence there are theoretical guarantees for these neural network architectures. 
A preliminary implementation of our proposed neural network architecture on FPGAs shows promising speed up compared to CPUs.
Beyond the numerical experiments in Section~\ref{sec:implementation}, we also tried some examples where the assumption (A2) is not satisfied. In these examples, we observed that our proposed neural network architectures also provide reasonable numerical outputs. These observations suggest that the assumption (A2) is sufficient but not necessary for our proposed architectures.

\begin{acknowledgements}
This research is supported by AFOSR MURI FA9550-20-1-0358. The authors also thank the Xilinx Center of Excellence at the University of Illinois, Urbana-Champaign UIUC to provide access to Xilinx Alveo boards and computing resources.
\end{acknowledgements}

\bibliographystyle{spmpsci}      
\bibliography{biblist}

\end{document}